\documentclass[a4paper,12pt]{article}

\usepackage{lipsum}
\usepackage{amsfonts}
\usepackage{amsmath}
\usepackage{bm}
\usepackage{graphicx}
\usepackage{epstopdf}
\usepackage{algorithmic}
\usepackage{fancyhdr}
\usepackage{hyperref}
\hypersetup{
    colorlinks=true,
    linkcolor=blue,
    filecolor=blue,      
    urlcolor=blue,
    citecolor=blue
}
\fancyheadoffset[RE,LO]{-0\textwidth}
\usepackage[right=3cm,left=3cm,top=3cm,bottom=3cm,headsep=1.5cm,footskip=1.5cm]{geometry}

\title{How do degenerate mobilities determine singularity formation in Cahn-Hilliard equations?
\thanks{Preprint submitted to arXiv.org on the 13th of January of 2021.}}

\author{Catalina Pesce\thanks{Mathematical Institute, University of Oxford, Woodstock Road, Oxford, OX2 6GG, UK. Email: pesce@maths.ox.ac.uk, 
muench@maths.ox.ac.uk}
\hspace{0.05cm} and 
Andreas~M\"unch\footnotemark[2]
}

\date{}

\usepackage{amsopn}

\usepackage{multirow}
\usepackage{booktabs}
\usepackage[labelfont=bf,font=scriptsize]{caption}
\usepackage[font=scriptsize]{subcaption}
\usepackage[makeroom]{cancel}
\usepackage{afterpage}
\usepackage[abs]{overpic}
\usepackage[colorinlistoftodos]{todonotes}
\usepackage{bm}

\counterwithin{equation}{section}
\counterwithin{figure}{section}

\newcommand{\eps}{\varepsilon}
\newcommand{\order}[1]{{O}(#1)}
\newcommand{\uinit}{u_{\text{init}}}

\newcommand{\rstar}{r_*}
\newcommand{\R}{\mathbb{R}}
\newcommand{\pphi}{\varphi}

\definecolor{BrickRed}{cmyk}{0,0.89,0.94,0.28}

\newsavebox{\gbox}
\newlength{\mylen}

\usepackage{comment}

\begin{document}

\maketitle

\begin{abstract}
Cahn-Hilliard models are central for describing the evolution of interfaces in
phase separation processes and free boundary problems.  In general, they have
non-constant and often degenerate mobilities.  However, in the latter case, the
spontaneous appearance of points of vanishing mobility and their impact on the
solution are not well understood.  In this paper we develop a singular
perturbation theory 
to identify a range of degeneracies for which the solution of the Cahn-Hilliard
equation forms a singularity in infinite time.  This analysis forms the basis
for a rigorous sharp interface theory and enables the systematic development of
robust numerical methods for this family of model equations.
\end{abstract}



\normalsize

\section{Introduction}
Since its introduction
\cite{cahn_free_1958,cahn_spinodal_1961,cahn_phase_1965,cahn_spinodal_1971},
the Cahn-Hilliard equation and its many variations have become fundamental
tools for describing the separation of phases over a large range of time and space
scales in many applications. In the most basic case of two
partially miscible materials, such as binary alloys or polymeric liquids,  this
includes the early onset of the phase separation from a homogeneous, unstable
state via spinodal decomposition \cite{cahn_spinodal_1961,cahn_phase_1965}, and
subsequent nonlinear evolution at later stages where coarsening occurs
\cite{gawlinski_domain_1989,rogers_numerical_1988} until
the quasi-stationary stages where only few, large and almost homogeneous domains remain.  Due also to its ability to allow for topological changes of the
domain, phase-field models based on extensions of the Cahn-Hilliard equation
are frequently used as the basis for numerical simulations of, for example, the
evolution of interfaces between immiscible liquids. These applications exploit
the fact that in a phase-field model, the interfaces are represented by a thin
layer over which the order parameter varies rapidly but continuously. 
Examples of such processes are surface diffusion and electromigration in
crystals and alloys \cite{cahn_surface_1994, taylor_linking_1994,
cahn_cahn-hilliard_1996, dziwnik_anisotropic_2017, bhate_diffuse_2000,
barrett_phase_2007, barrett_finite_2004}, motion of immiscible fluids with free
boundaries \cite{ding_diffuse_2007, abels_thermodynamically_2012,
hosseini_isogeometric_2017,  boyer_study_2006,sibley_unifying_2013}, polymer
blends \cite{mcguire_kinetics_1995,de_gennes_dynamics_1980,
castellano_mechanism_1995},  tumour growth models
\cite{cristini_nonlinear_2009,oden_toward_2015,lima_selection_2017} or
lithiation in battery electrodes  \cite{meca_localized_2018}, to name just a
few.  

Stated in the form introduced by \cite{cahn_free_1958,cahn_spinodal_1971}, 
the Cahn-Hilliard equation can be written as 
\begin{subequations}\label{eqn:chgen}
\begin{equation}\label{eqn:chgena}
u_t=-\nabla\cdot\mathbf{j}, \qquad
\mathbf{j}=-M(u)\nabla\mu, \qquad
\mu=-\eps^2\nabla^2 u + f'(u),
\end{equation}
with the (conserved) order parameter $u$, such that $|u|\leq 1$ and $\eps>0$.
The homogeneous free energy and mobility are, respectively
\begin{equation}
\label{eqn:logfe-mobdeg1}
f(u)=(\theta/2)\left[(1-u)\ln(1-u)+(1+u)\ln(1+u)\right]+(1-u^2)/2,\qquad
M(u)=1-u^2,
\end{equation}
\end{subequations}
where $\theta\geq 0$ denotes a normalised temperature. For $0<\theta < 1$ the free energy in \eqref{eqn:logfe-mobdeg1} has two
distinct minima $u_\pm$ and the system separates into two phases with those relative concentration values. 
The diffuse interface layers between these phases domains are thin if $\eps$ is small.
The term Cahn-Hilliard is often used more broadly to describe a class 
of phase-field models that have
the general form \eqref{eqn:chgen} but different free energies and mobilities,
for example, a quartic polynomial with fixed minima $u=\pm 1$, 
such as
\begin{equation}\label{qfe}
f(u)=\frac 12 (1-u^2)^2, 
\end{equation}
and a constant mobility $M(u) \equiv 1$. 

While a different mobility does not change the energy landscape, it does strongly affect
the kinetics of the process. For constant mobility, the flux depends only on the gradient
of the chemical potential $\mu$ and the diffusive flux $\mathbf{j}$ can freely transport 
material through the bulk in the direction of decreasing $\mu$. The kinetics become clearer when one takes $\eps$, and thus the interface width, to $0$. For the constant mobility Cahn-Hilliard equation, Pego \cite{pego_front_1989} showed, via matched asymptotics, that the sharp interface limit is the Mullins-Sekerka problem, which inspired the rigorous proof by Alikakos, Bates and Chen \cite{alikakos_convergence_1994}. The Mullins-Sekerka problem couples the interface motion to the bulk diffusion between the domains at the late stages of the coarsening process. 

In contrast, nonlinear mobilities that degenerate at or near the minima of the free energy suppress
bulk diffusion, so that transport along the interface, i.e.\ surface diffusion, becomes more
important. Using asymptotic methods, Cahn et. al \cite{cahn_surface_1994} demonstrated that for $\theta\ll 1$ and for a double-obstacle free energy ($\theta=0$), the sharp interface limit 
for \eqref{eqn:chgen} is simply the surface diffusion equation, with no
transport across the bulk, at least to leading order. On the other hand, for the case
of a quartic free energy (\eqref{eqn:chgena} and \eqref{qfe}), 
the degenerate mobility leads to a subtle balance between bulk and
surface diffusion, so that to leading order, both enter the sharp interface
limit \cite{dai_coarsening_2014, lee_sharp-interface_2016}. This has come as a
surprise to some in the community, as by  routine application of Pego's
asymptotic approach, one can easily miss the contribution from bulk diffusion
\cite{lee_degenerate_2015, lee_response_2016, voigt_comment_2016} and obtain
the wrong sharp interface model.  The correct and consistent evaluation of the
flux requires the use of exponential matching
\cite{lee_sharp-interface_2016}.

The results in 
\cite{lee_degenerate_2015,lee_sharp-interface_2016,lee_response_2016} highlight 
the subtleties arising from degenerate mobilities and
the importance of investigating the equations carefully. 
Besides the derivation of the sharp interface limit in \cite{lee_sharp-interface_2016}, 
another aspect became apparent upon 
solving the axisymmetric PDE on
a circular, two-dimensional domain with initial data $u_0$ strictly bounded
between $-1$ and $+1$. This represents a situation where the phases have separated
into two domains, a disc centered at the origin with a composition close to one
phase, surrounded by an annular region with a composition near
the opposite phase, and a diffuse interface between them. 
As a numerical result in \cite{lee_sharp-interface_2016} reveals, the solution
evolves so that near the interface, $|u|$ develops a maximum that quickly
approaches $1$, that is, the value for which the mobility degenerates. 

This phenomenon is intimately connected with a property of stationary solutions
to Cahn-Hilliard equations with smooth polynomial free energies.
These are well-studied and, in particular, existence and uniqueness of 
small energy stationary solutions has been proven in 
Niethammer~\cite{niethammer_existence_1995} using a rigorous matched asymptotics expansion technique that also
captures the qualitative features of the solution. Nontheless the following property (referred to
as the Gibbs-Thompson effect in \cite{dai_coarsening_2014}) has often been
overlooked: In the presence of curved interface layers between 
phases, the chemical potential is non-zero in equilibrium, and the 
``outer'' solution i.e.\ the solution away from interface layers, 
differs from the minima of the free energy by a small
amount proportional to the curvature. Inside convex domains, the value is in
fact outside of the interval delineated by the minima of the free energy (i.e.\
here $\pm 1$). Since time dependent solutions of Cahn-Hilliard equations
monotonically decrease their energy, they are expected to converge to
stationary solutions; in particular, to the one investigated by
Niethammer~\cite{niethammer_existence_1995}. As a result, $u$ must approach
$\pm 1$ somewhere, thus forcing the degenerate mobility $M(u)$ 
in \eqref{eqn:logfe-mobdeg1} to become zero.

This observation raises interesting questions that have important
implications for established practices. To begin
with, is $u=\pm 1$ achieved in finite or infinite time? What determines this? Since at those points mobilities like \eqref{eqn:logfe-mobdeg1} degenerate, does it depend on how degenerate $M$ is? 
We note that some authors \cite{mahadevan_electro_1999, ratz_surface_2006, torabi_inter_2012, torabi_new_2009, gemmert_phase_2005, wolterink_spinodal_2006} choose a low-degeneracy mobility with a degree of two, that is, the square of the form used in \eqref{eqn:logfe-mobdeg1}, but higher degeneracies can also be useful to understand the full spectrum of the solution's behaviour. Next, what happens for example in the case that $|u|$ approaches 1 (and hence $M(u)$
approaches 0) in finite time? Can this lead to loss of regularity and, thus, to singularity formation? Will the vanishing mobility freeze the solution
there and prevent the set $\{M(u)=0\}$ to move, and how will that influence the
evolution of the diffuse interface and hence the sharp interface limit?  How
will that affect long-time pattern formation in numerical simulations? 

An early
paper by Elliott and Garcke \cite{elliott_cahnhilliard_1996}, where they prove existence
of solutions for a class of degenerate Cahn-Hilliard models, first raises the
question of how the set $\{M(u)=0\}$ evolves. In fact, pinning was observed in numerical solutions of degenerate Allen-Cahn/Cahn-Hilliard systems \cite{barrett_finite_2002}, in contrast to the constant mobility case. Moreover,
solutions with a waiting time behaviour are also conceivable \cite{novick-cohen_thin_2011}.
Numerical experiments in \cite{barrett_finite_1999} demonstrate that the choice
of the relative magnitudes of the mesh and the temporal step size yields at
least two solutions with very different behaviour. If the mesh sizes are taken
to zero much faster than the step size, the solution the scheme converges to is
pinned at the boundary of the set $\{M(u)=0\}$ and hence it is stationary, while
another, moving, solution emerges if the step and mesh sizes are
in a distinguished limit with each other.  Such a behaviour is important to
know and understand, as selecting and changing step and mesh sizes is standard
practice in numerical simulations, and in fact is often done automatically as
part of the adaptivity implemented in ready-to-run simulation packages. These results highlight the role of non-uniqueness of solutions, which also prompts questions
about the implications for the sharp interface limit, in particular, can different solutions
have different sharp interface limits? As a consequence, it becomes essential to investigate
the situation at points $\{M(u)=0\}$, which are typically points where the solution becomes singular, in the sense that the regularity is reduced
\cite{constantin_droplet_1993}, requiring the introduction of weak
solution concepts \cite{elliott_cahnhilliard_1996,bernis_higher_1990}.

These questions overlap with another important 
class of fourth order PDEs with degenerate
mobility, namely the surface tension driven thin film model. 
The thin film equation for surface tension driven flows (and its variants) has a rich
history and a huge literature covering questions (1) on the formation and
evolution of sets where the mobility is zero  \cite{constantin_droplet_1993,
constantin_singularity_2018, bertozzi_singularities_1994, bernis_higher_1990,
bertozzi_lubrication_1994, bertozzi_symmetric_1996,
vaynblat_rupture_2001,vaynblat_symmetry_2001, witelski_dynamics_2000,
erneux_nonlinear_1993, peschka_thin-film_2010, peschka_self-similar_2010,
munch_numerical_1999},
especially in the context of the fundamental questions in fluid mechanics
about the moving 
three-phase contact-line, see e.g.\ \cite{greenspan_motion_1978,
barenblatt_problem_1997, lacey_motion_1982,hocking_spreading_1983}, see also
\cite{bertozzi_mathematics_1998}, and (2) on the impact of the degree of
degeneracy on the solutions \cite{almgren_stable_1996,king_moving_2001}.
In a situation where sufficient pressure is applied to the boundary of a thin film, singularities are always forced. In particular, questions about how the
thin film height goes to zero and if the singularity occurs in finite or infinite time have been discussed in the literature \cite{constantin_droplet_1993,
constantin_singularity_2018, bertozzi_singularities_1994}.  The connection
between thin films and degenerate Cahn-Hilliard problems offers a lot of
potential, and even though it has been highlighted earlier
 \cite{thompson_emerging_2000, novick-cohen_thin_2011}, it has rarely been directly
exploited. Furthermore, in this paper, the asymptotic techniques from the
thin film literature have proven fundamental for understanding the behaviour of
a degenerate Cahn-Hilliard problem.

The answer to the questions above depend on the details of the problem.  While
the originally physically motivated problem suggests a logarithmic free energy
or its deep quench limit, that is the obstacle free energy
\cite{cahn_cahn-hilliard_1996, cahn_spinodal_1961}, which remains the
preferred choice for part of the community in the context of surface diffusion,
many others have chosen polynomial double well potentials in combination with
degenerate mobility of up to second degree degeneracy instead.  
Typically, in these cases, $M$ vanishes exactly at the minima of the double
well potential, while in \eqref{eqn:logfe-mobdeg1}, the minima of $f$ are 
not located at $u=\pm 1$ for $0<\theta<1$.
The choice of a quartic polynomial is also frequent where Cahn-Hilliard
equations are coupled with the Navier-Stokes equation to solve free boundary problems
numerically, and also in one-sided models where significant bulk diffusion
is only desired on one side of the interface \cite{glasner_diffuse_2003,
lu_diffuse-interface_2007, dai_motion_2012,dai_weak_2015}.  

In this paper,  we
continue with this class of problems and focus on a particular domain, namely we consider equation \eqref{eqn:chgena} on a
two-dimensional axially symmetric domain with a quartic polynomial free
energy $f$ and a mobility $M$ which vanishes at the minima of $f$, and where the
degree of degeneracy is treated as a parameter $n \in \R_+$.  
We show that the solution develops points
where $|u|\to 1$, which can form in either finite or infinite time, and that for a range of
mobilities, there are attracting solutions that belong to the latter category.
These solutions are analysed by singular perturbation methods.

The layout of the paper is as follows. In section 2, we summarise the precise
statement of the axially symmetric Cahn-Hilliard equation that we consider. In
section 3, we present the result of numerical solutions for a range of values
$n$ for the degree of degeneracy of $M$.  In particular, we determine the
self-similar regions that develop in the long-time solution for the example of
$n=4$. In section 4, we systematically derive an
asymptotic approximation for the long-time behaviour of this solution for the
case $n>2$ using matched asymptotics, giving, in particular, the similarity
exponents for the different asymptotic regimes of the solution.  In section 5,
we discuss our results and point to possible further questions and avenues of
research.

\section{Problem statement} We take \eqref{eqn:chgena} on the 2-dimensional unitary ball for a radially symmetric smooth function $u=u(r,t)$, which written in polar coordinates corresponds to
\begin{subequations}\label{tp}
\begin{equation}\label{tppde}
u_t=-\frac1r\frac{\partial(r j)}{\partial r},
\qquad
j=-M(u)\frac{\partial \mu}{\partial r},
\qquad
\mu=-\frac{\eps^2}{r}\frac{\partial}{\partial r} \left(r\frac{\partial u}{\partial r}\right)
+ f'(u),
\end{equation}
for $0<r<1$ and $t>0$. 
We also assume that $u$ and its derivatives with respect to Cartesian coordinates are continuous
at the origin, which implies the boundary conditions
\begin{equation}\label{tpbcat0}
\partial_r u =0,\qquad \partial_{rrr} u=0, \qquad \text{at } r=0.
\end{equation}
and moreover, we assume that we have a neutral surface at $r=1$ and no flux, so that
\begin{equation}\label{tpbcat1}
\partial_r u=0, \quad j=0, \qquad \text{at } r=1.
\end{equation}
We also need to prescribe an initial condition
\begin{equation}\label{ui}
u(r,0)=\uinit(r),
\end{equation}
which we specify further in the next section. Typically, it will be
a scaled tanh-profile that is strictly bounded between $+1$ and $-1$.
The homogeneous free energy is given by a double-well potential 
\begin{equation}\label{fes}
f(u)=\frac 12 (1-u^2)^2, 
\end{equation}
and the mobility by
\begin{equation}\label{mobg}
M(u)=(1-u^2)^n,
\end{equation}
\end{subequations}
where $n> 0$ is a real parameter.
This is the form of the mobility that we shall use most often in this paper,
though we also discuss two variants, in particular where we refer to results
in the literature. One variant, considered
for example by Elliott and Garcke \cite{elliott_cahnhilliard_1996}, is to
truncate the mobility, so that
\begin{equation}\label{mobp}
M(u)=(1-u^2)_+^n
\end{equation}
where the subscript $+$ denotes the positive part of the expression (taken before
the expression is raised to the power of $n$). A second variant, used by 
Dai and Du \cite{dai_coarsening_2014, dai_weak_2015, dai_computational_2016}, 
is to take instead the absolute value,
\begin{equation}\label{mobabs}
M(u)=|1-u^2|^n.
\end{equation}
In all cases, the parameter $n>0$ determines the degree of degeneracy of the mobility. These variants differ for values of $u$ with $|u|>1$, where
\eqref{mobg} becomes undefined for non-integer $n$, or changes sign for odd
$n$, leading to ill-posedness unless $n$ is even. We avoid these situations in
the current article by focusing on bounded solutions $|u|\leq 1$.

\section{Numerical solution}\label{numerical-results}

We begin by inspecting numerical solutions of \eqref{tp} for three different
groups of the mobility: (a) constant mobility, $M=1$, corresponding to
\eqref{mobg} with $n=0$; (b) degenerate mobility with $n=1$; (c) degenerate
mobility with $n=4$. The initial data is given by \eqref{ui}, with 
\begin{equation}\label{uinit}
\uinit(r)=-0.95\tanh \left(\frac{r-0.5}{\eps}\right).
\end{equation}
Such tanh-like profiles are a common choice for numerical simulations where 
phase-field models are used to track the evolution of a free interface, or to
capture the late stages of a phase separation process, see \cite{torabi_new_2009}, \cite{torabi_inter_2012}.
Unless otherwise stated, we choose
\begin{equation}
\eps=0.1.
\end{equation}
Notice the initial profile above satisfies the boundary conditions \eqref{tpbcat0},
\eqref{tpbcat1} except for exponentially small terms. Replacing \eqref{uinit}
in the vicinity of $r=0$ and $r=1$ by constant values $+1$ and $-1$, respectively, so that the initial profile satisfies the boundary conditions, did not change the 
numerical results in any noticeable way.

The numerical solutions for \eqref{tp} with initial data \eqref{uinit} presented here were obtained
by a finite difference code using centred differences in space and
and implicit Euler scheme in time. The spatial grid was equidistant, and we used
a step doubling scheme to control the error in time.

\begin{figure}
\savebox{\gbox}{\includegraphics[clip=true,width=0.52\textwidth]{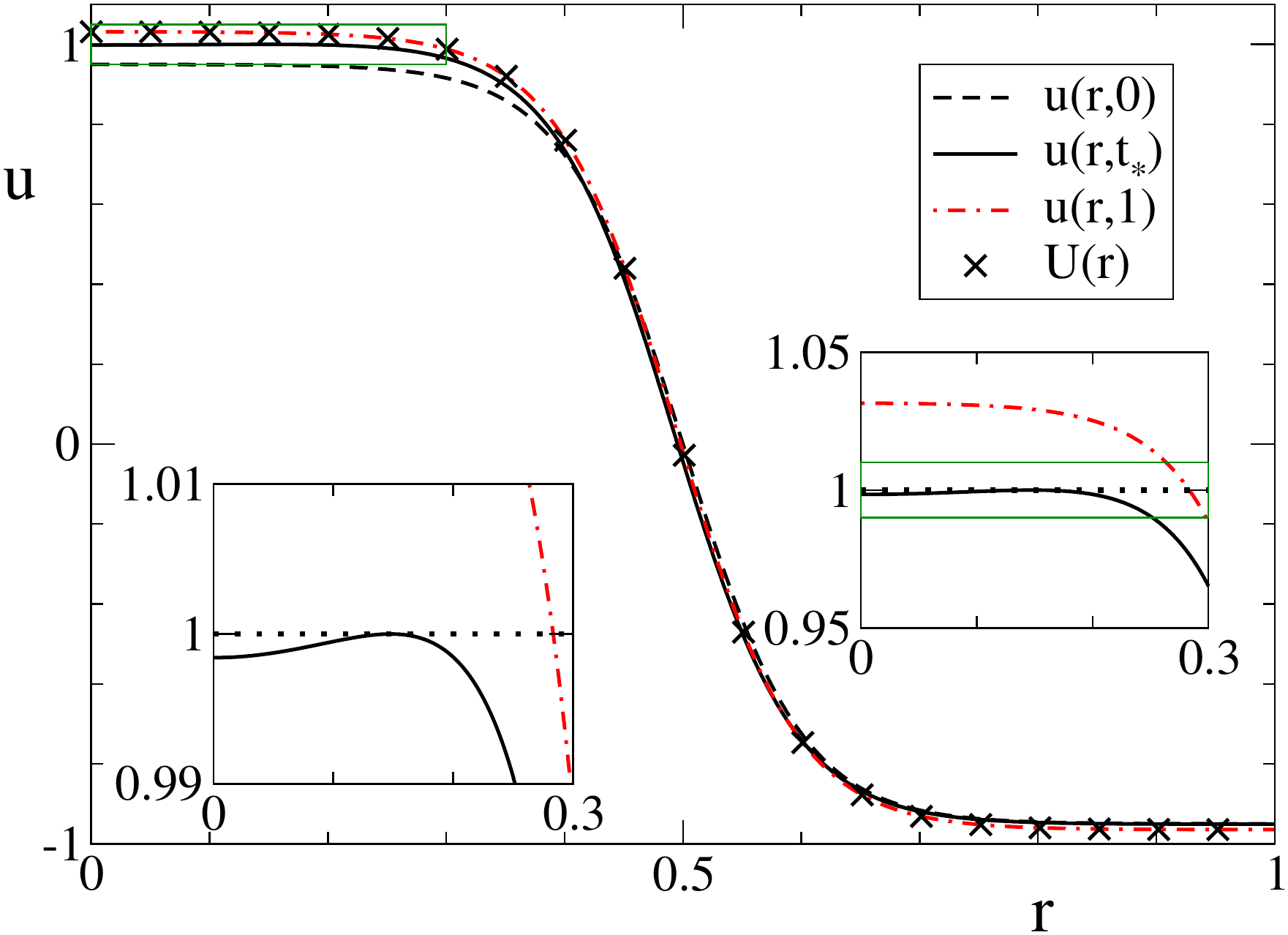}}
\settoheight{\mylen}{\usebox{\gbox}}
\hfill%
\usebox{\gbox}%
\hfill%
\parbox[b]{0.35\textwidth}{%
\includegraphics[clip=true,width=0.35\textwidth,height=0.5\mylen]{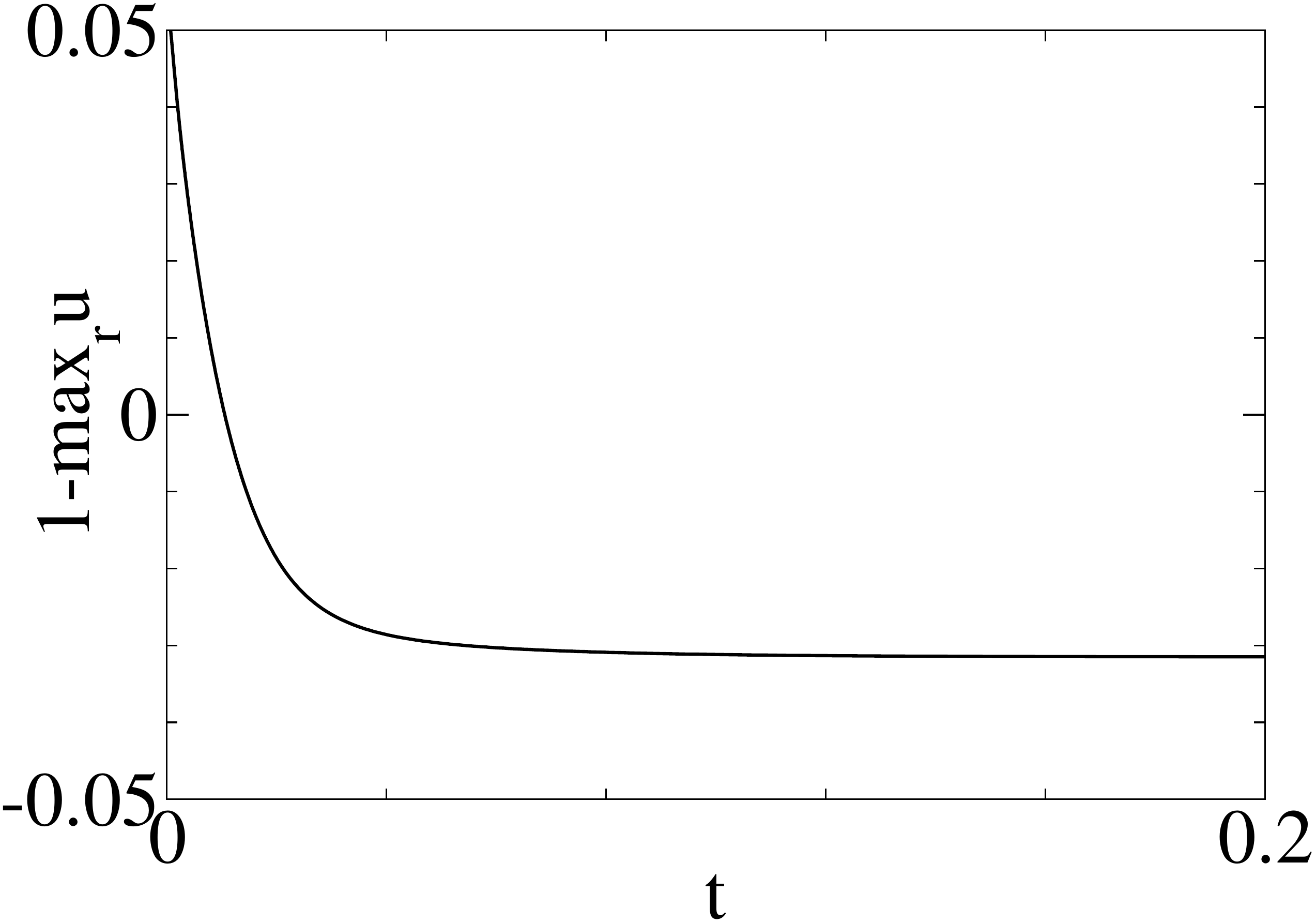}\\
\includegraphics[clip=true,width=0.35\textwidth,height=0.5\mylen]{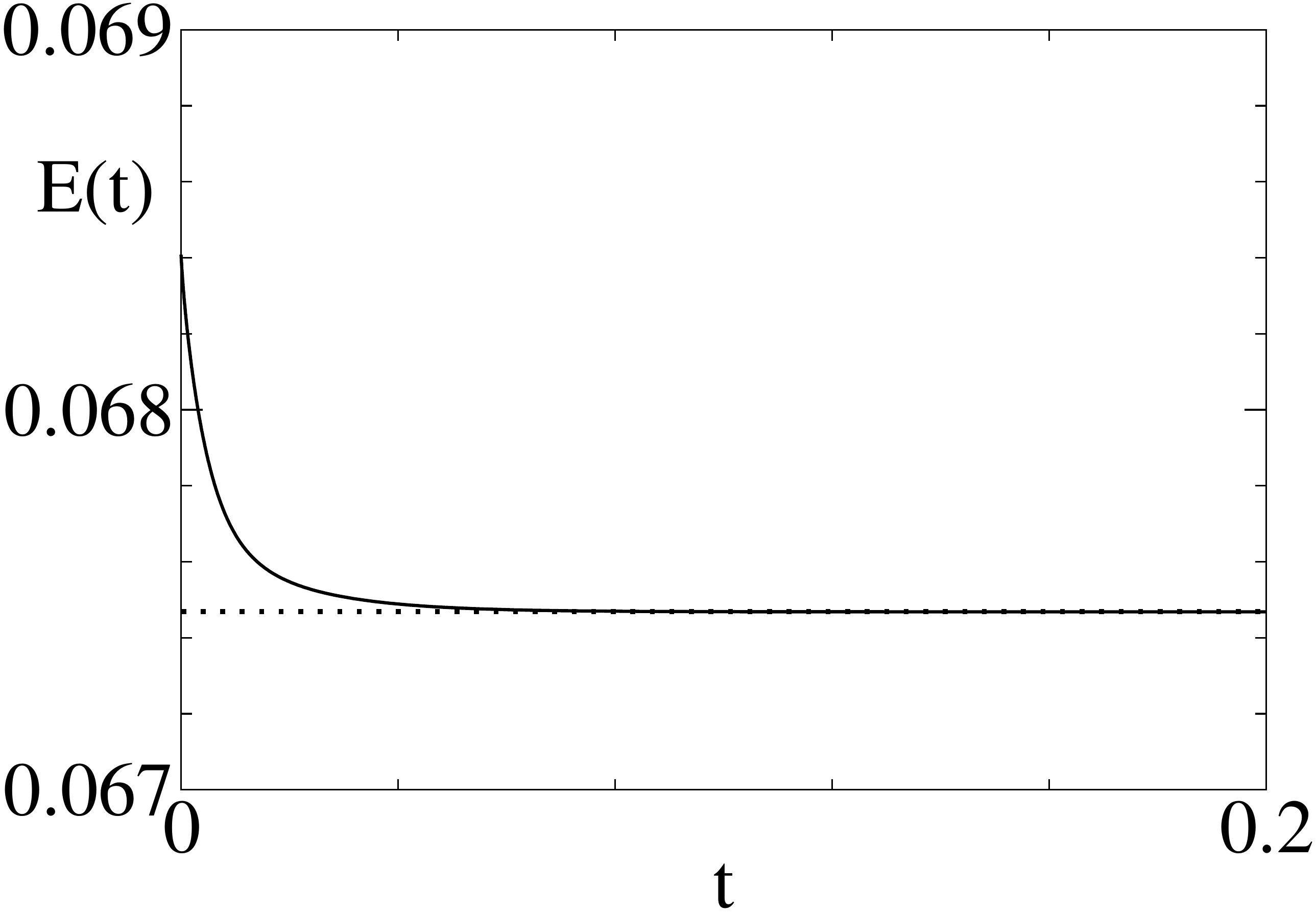}}
\hfill%
\caption{\label{fig:constmob} (a) Left: Solution $u$ of \eqref{tp} with
constant mobility $M=1$, at different times $t=0$, $t=t_*=1.06\times 10^{-2}$,
$t=1$. This is compared with the stationary solution $U(r)$, which satisfies
\eqref{stat} and \eqref{mU}. The right inset shows a zoom of the area
delineated by a thin solid box, and the left inset in turn zooms into the area
between the horizontal vertical lines in the first inset. A thin dotted line at
$u=1$ has been added to both insets for guidance. (b) Right top: Evolution of
$1-\max_r u(r,t)$ and bottom: of the energy $E(t)$ of the solution. The thin
dotted line shows the energy for the stationary solution of
\eqref{stat}, \eqref{mU}.}
\end{figure}

\paragraph{Constant mobility}
We see in fig.~\ref{fig:constmob}(a) that the solution develops 
a maximum at a radius $\bar r(t)$ near $r=0.15$,
which quickly crosses $u=1$ at $t=t_*=0.0106$, after which $u$ settles into a
stationary solution.  
To understand better the intuition behind the long time solution of the constant mobility case we need to introduce the free energy associated with this system, which, in polar coordinates, is given by
\begin{equation}\label{fed}
{\cal F}[u]=\int_0^1 \left[\frac{\eps^2}2 (\partial_r u)^2+f(u) \right] r\mathrm dr.
\end{equation}
This energy is always non-increasing along a solution trajectory, i.e., if $u(t)$
is a solution of \eqref{tppde}, \eqref{tpbcat1}, then $E(t):={\cal F}[u(t)]$ satisfies
\begin{equation}
\frac{\mathrm dE}{\mathrm d t}= -\int_0^1 M(u) \left(\partial_r \mu\right)^2 r\mathrm dr \leq 0.
\end{equation}
Since $E\geq 0$, this means that $E\to E_\infty$ as $t\to \infty$, and $
{\mathrm dE}/{\mathrm d t}\to 0$. 
(Notice this is true also for general nonlinear mobilities provided they are non-negative).
Since $M$ is constant, $\mu$ converges to a constant, say
$\mu_c$. If the solution converges to a stationary solution $U(r)$ of
\eqref{tp}, which is known to be true at least for the case of constant mobility with logarithmic \cite{abels_equilibrium_2007} or quartic polynomial homogeneous free energy \cite{rybka_convergence_1999},
then $U(r)$ must satisfy
\begin{subequations}\label{stat}
\begin{align}
-\frac{\eps^2}r\frac{\mathrm d}{\mathrm d r}\left(r\frac{\mathrm d  U}{\mathrm d r}\right)
+f'(U)&=\mu_c,\\
U'(1)&=0,\\
U'(0)&=0.
\end{align} 
\end{subequations}
The additional degree of freedom $\mu_c$ is used to accommodate a mass constraint that the
solution inherits from the initial condition. The system \eqref{tp}, \eqref{tpbcat0}, \eqref{tpbcat1}
conserves mass, that is, for
\begin{equation}
m(t):=\int_0^1 u(r,t) r\mathrm dr
\end{equation}
one easily finds that $\mathrm d m / \mathrm d t=0$ along a solution $u(r,t)$, therefore, 
for the stationary solution, we need to enforce 
\begin{equation}\label{mU}
\int_0^1 U(r) r\mathrm d r=m_0,
\end{equation}
where 
\begin{equation}\label{mzero}
m_0= \int_0^1 u_{\text{init}}(r)rdr.
\end{equation}
It has been shown in \cite{niethammer_existence_1995}
that for any initial mass $m_0$,
there exists a unique pair of solutions to
\eqref{stat}, in the set of smooth functions with sufficiently small energy ${\cal F}[U]=\order\eps$, 
which are identical up to a reflection $u\rightarrow -u$. Hence for initial data with small enough energy,
we expect the solution of the initial boundary value problem to converge to
one of these stationary states, namely the one closer to the initial data.
We can check this by superimposing the solution for \eqref{stat} onto
the long-time profile for the initial boundary value problem.
In addition, the stationary solution $U(r)$ exceeds $1$ by an $\order\eps$ amount.
This fact, which is a manifestation of the Gibbs-Thomson effect
\cite{dai_coarsening_2014, lee_sharp-interface_2016} 
is often missed, but it is an important feature of the evolution.  It basically precludes the
possibility that $|u|<1-\delta$ point-wise, for some positive $\delta$. If that were the case, the small-energy stationary profile to which the
solution converges would be strictly bounded by $|u|<1$, contradicting
the result in \cite{niethammer_existence_1995}.  This is a strong indication that 
$|u|\to 1$ at some point(s) $r_*$ and at either a time $t_* <\infty$ or at the limit $t\to\infty$. 
We refer to these two cases as finite-time and infinite-time touchdown, respectively.
For constant mobility, the former of the two occurs.

\paragraph{Degenerate mobilities with $n=1$}

\begin{figure}
\includegraphics[clip=true,width=0.45\textwidth]{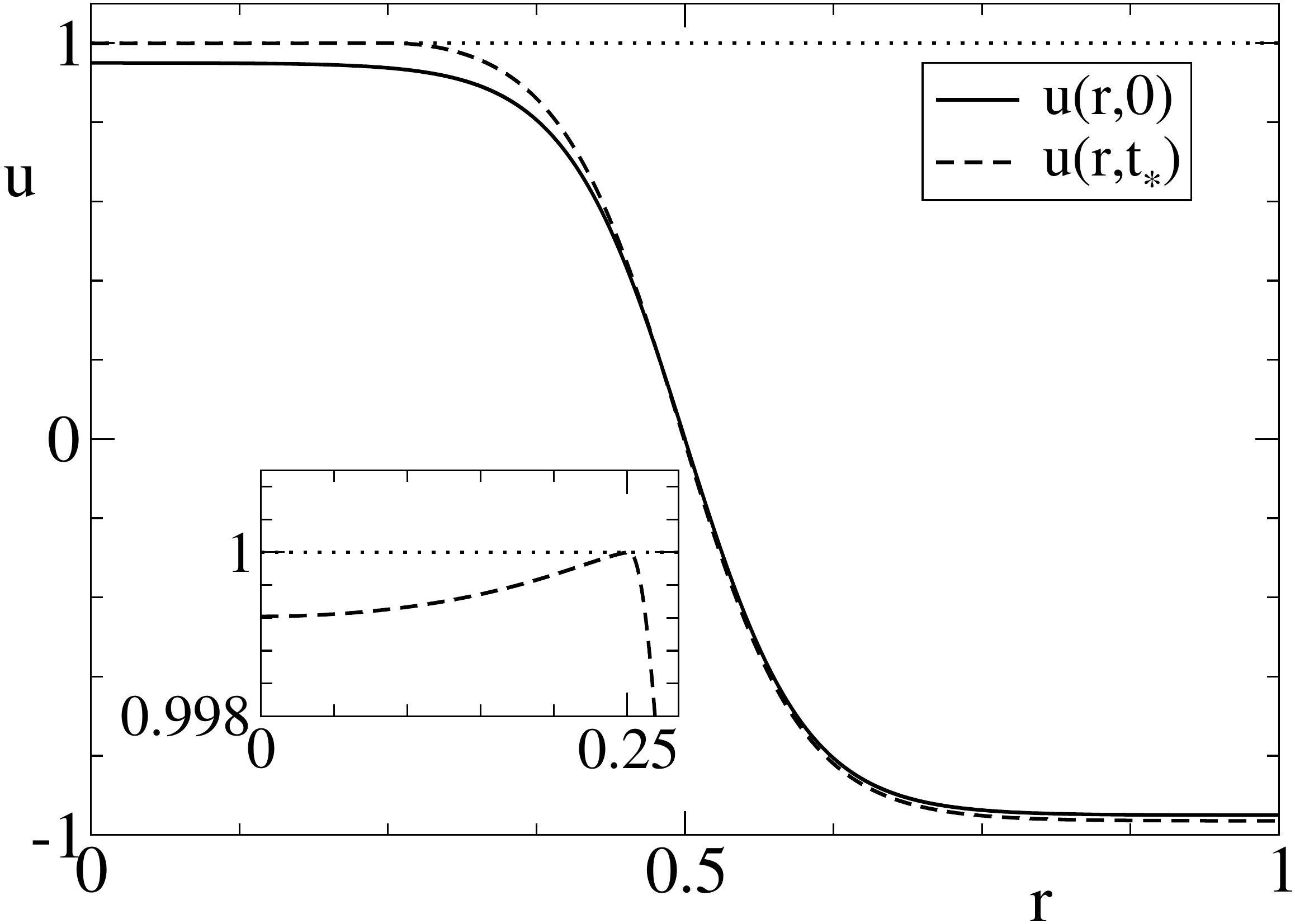}
\qquad
\includegraphics[clip=true,width=0.45\textwidth]{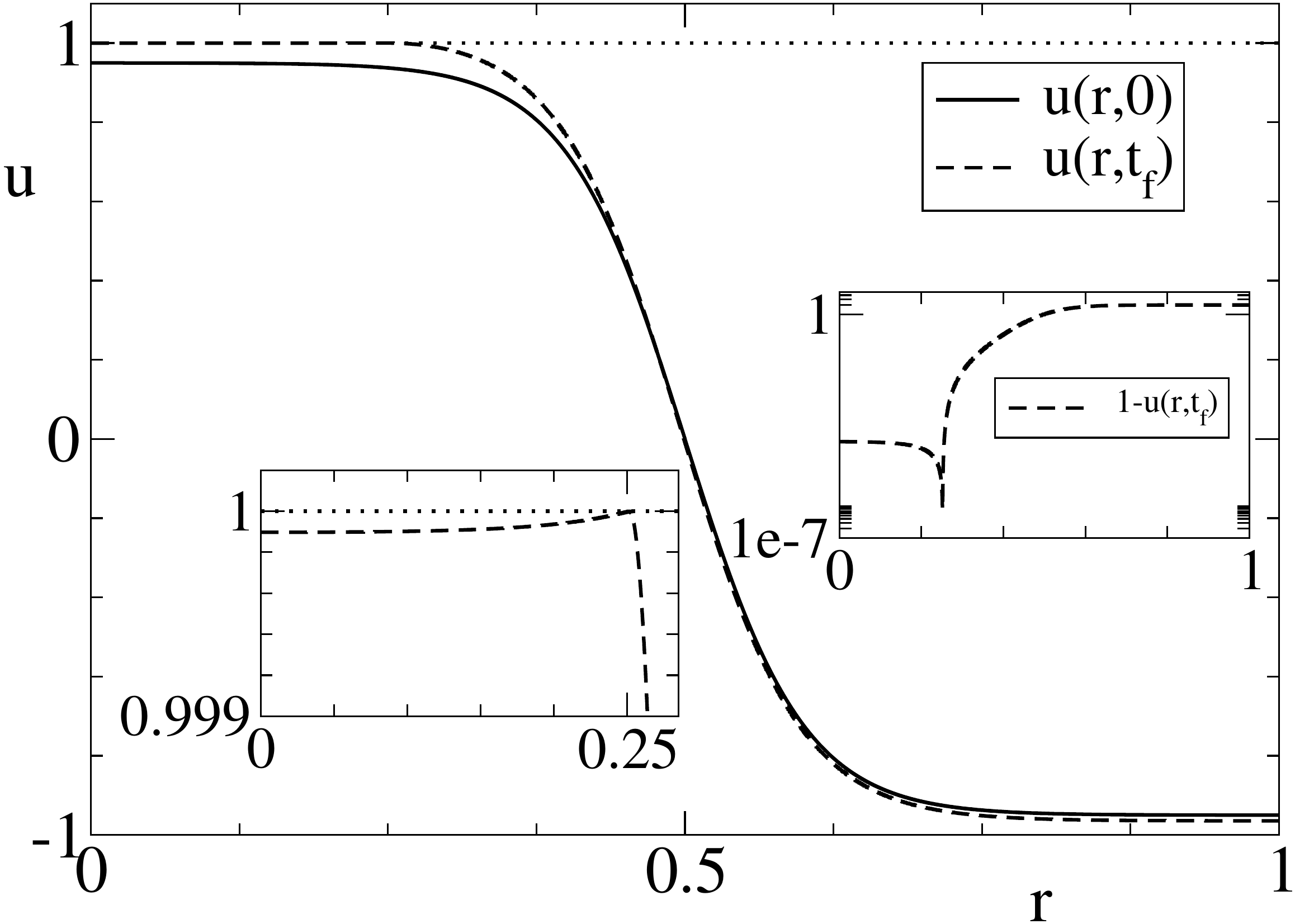}
\caption{\label{fig:degmobn1} 
(a) Left: Solution $u$ of \eqref{tp} with degenerate mobility $n=1$,
for the initial time ($t=0$) and after the touchdown, 
as illustrated by the zoom in the inset.
\eqref{mU}. (b) Right: Solution of \eqref{tp} with mobility $n=4$,
for the initial time and for $t=t_{final}=10^{15}$, when we stopped
the simulation (though it could still have been continued). The maximum of
$u$ is very close to the $u=1$, see bottom left inset, but, as a semi-log
plot of $1-u$ in the top right inset reveals, has touched there.}
\end{figure}

The numerical solution behaves similarly to
the constant mobility case in that $\max_r u$ reaches $1$ in finite time, at
$t=t_*=3.44$ 
for a grid spacing of $\Delta r=10^{-4}$, 
see fig.~\ref{fig:degmobn1}(a) and  the inset. 
However, we noticed that $t_*$ increases significantly upon decreasing
$\Delta r$, leaving open the possibility that this finite-time touchdown is a
numerical artefact. This means that we could still be approximating a solution that
stays away from $|u|=1$ for all times, and are only thrown off this trajectory
due to numerical errors. 

Finite time touchdowns - or in fact crossings -- of $|u|=1$ have been
reported in the literature. 
For example, for the fully two-dimensional simulations in \cite{dai_computational_2016},
where the authors use an absolute value mobility \eqref{mobabs}, the solution $u$ crosses the
bound $|u|=1$ on the convex
side of interfaces between different phases and converges to a quasi-stationary
profile with a larger value than the one as predicted by the Gibbs-Thomson effect.
Other analytical results \cite{elliott_cahnhilliard_1996} 
prove existence of solutions $u\leq 1$, which allow
for touchdowns but not for crossings into $|u|>1$.
The occurrence of qualitative different behaviours for one initial value
problem is consistent with the non-uniqueness of solutions that is known
to occur for initial value problems of high-order degenerate parabolic PDEs
in the Cahn-Hilliard and thin-film context, see for example
\cite{dai_weak_2015, dai_motion_2012, king_moving_2001, beretta_nonnegative_1995, 
dal_passo_existence_1999}. 

Moreover, where the solution
achieves a value for which the mobility vanishes, it typically loses
regularity.  This can be illustrated by a formal argument made for a Hele-Shaw
model in \cite{constantin_droplet_1993} (see also a rigorous version in
\cite{constantin_singularity_2018} and a similar argument for a more general thin film problem in \cite{bertozzi_singularities_1994}), applied here to $v:=1-u$. Let $n=1$, at the minimum
$r=\bar r(t)$ of $v$ at time $t$, which we denote $\bar{v}(t):=v(\bar{r}(t),t)$, we have that
\begin{equation}
\frac12\frac{\mathrm d}{\mathrm d t}
\ln\left(\frac{2-\bar{v}}{\bar{v}}\right)=\left.\frac1r\partial_r\left(r\partial_r \mu\right)\right|_{r=\bar{r}(t)}.
\end{equation}
When $\bar{v}(t)\to0$, the left hand side blows up and the second derivative of $\mu$, and hence the fourth derivative of
$v$ (or $u$), must do so too. The argument can also be extended to $n>1$ by computing instead $\frac{\mathrm d}{\mathrm d t}
\left(\bar{v}^{1-n}\right)$.  As a
consequence, rigorous treatments of \eqref{eqn:chgena}, \eqref{fes} and
\eqref{mobp} or \eqref{mobabs} introduce weak formulations and then typically
prove existence of such solutions via regularisation of the degeneracy, with
different outcomes depending on the details of the weak formulation and the
regularisation method.  In
\cite{elliott_cahnhilliard_1996} for example, the authors prove existence of
solutions for $n\geq 1$ that satisfy $|u|\leq 1$, using a regularised version
of \eqref{mobp}, while in \cite{dai_weak_2015}, Dai \& Du introduce a weak
solution concept and a regularisation of the mobility \eqref{mobabs} that
allows for solutions where $|u|$ can exceed 1. This is consistent
with the solutions they present in their numerical study 
\cite{dai_computational_2016}.

At this stage one may ask if the vanishing of the mobility along the solution
 can be avoided, so that, for example, the existence of
solutions with $|u|\leq 1$ by Elliott \& Garcke
\cite{elliott_cahnhilliard_1996} can be strengthened to show the existence of
a solution $u$ for which $|u|$ stays strictly
below 1 even in the limit as $t\to\infty$.
However, in the preceding section on the constant mobility case, we
gave an argument that rules out convergence to a stationary solution with
modulus less than $1-\delta$ for some $\delta>0$, which also carries over to the
degenerate case $n>0$. 
This implies that the solution (provided it converges to a stationary solution)
either achieves $\max_r |u|=1$ in finite time or, converges as $\max_r |u|\to
1$ in infinite time.

\paragraph{Degenerate mobilities with $n=4$}

As before, a maximum forms in the numerical solution that approaches $u=1$ at
some point $r=\bar r(t)$, but $1-u$ remains positive over many decades
of $t$. In fig.~\ref{fig:degmobn1}(b), $u$ still has not touched $u=1$ at
$t=10^{15}$. This suggests that the singularity is only approached in
infinite time. Moreover, the PDE
remains strictly parabolic and hence we expect it to have a unique classical solution,
that is, the same evolution should emerge for any other convergent numerical
scheme.  

In the following, we investigate the behaviour shown by this third example in
more detail numerically and via asymptotic analysis for the long time limit
$t\to\infty$, to conclude that the numerical solutions of \eqref{tp},
\eqref{mobg}, \eqref{uinit} converge to a leading order asymptotic
approximation that touches down in infinite time (and remains bounded away from
$|u|=1$ for any finite value of $t$.)

\subsection{Self-similar regions}

\begin{figure}
\includegraphics[clip=true,width=0.45\textwidth]{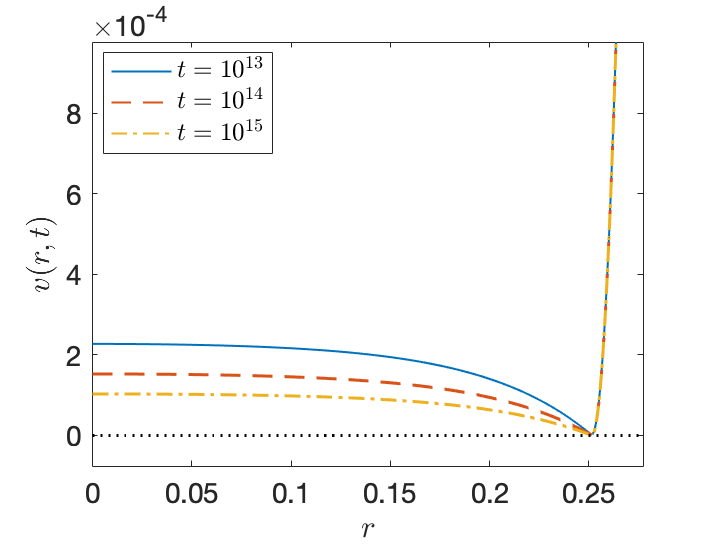}%
\qquad%
\includegraphics[clip=true,width=0.45\textwidth]{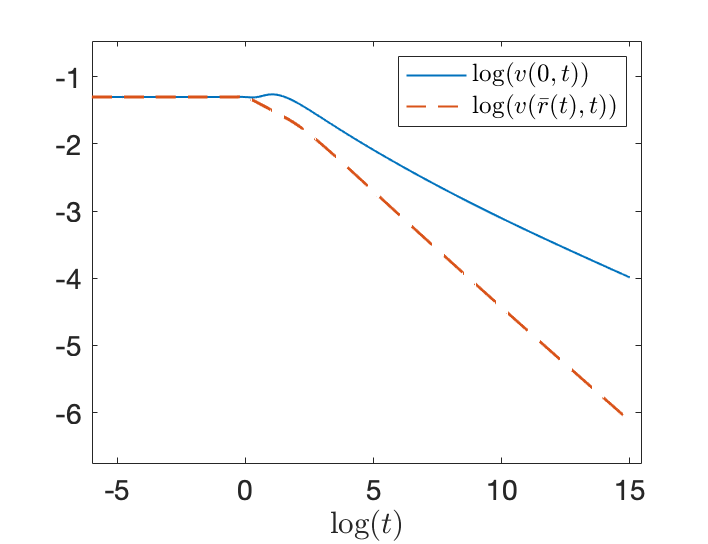}
\caption{\label{fig:zoomn4} (a) Left: The profile $v=1-u$ for the
solution of \eqref{tp} with mobility $n=4$,
for different times. The dotted line represents the $0$ value line.
(b) Right: Evolution 
of $v(0,t)$ and $v(\bar{r}(t),t)$ in a log-log plot.}
\end{figure}

We consider numerical results for three different values of $n=3$, 4, 5,
with the aim of investigating the structure of the solution
at large times.
To characterise the evolution as $t\to \infty$, we let $v=1-u$
and zoom into the regions of $r \in (0,1)$ where $|v|$ is small.  We first observe
that the region of $v$ closer to $r=0$, which we refer to as the central
region, evolves differently from the touchdown region near $r=\bar r(t)$.
There, the solution has a pronounced minimum $v(\bar r(t),t)$, and the
function decreases more rapidly than for $v(0,t)$.  In fact, the log-log plot
in fig.~\ref{fig:zoomn4} suggests that
$v(0,t)$ and $v(\bar r(t),t)$ display a power law behaviour for large $t$. Furthermore, both regions keep their qualitative shape, prompting us to seek self-similar
solutions with power-law scaling factors.

\begin{figure}
    \includegraphics[clip=true,width=0.45\textwidth]{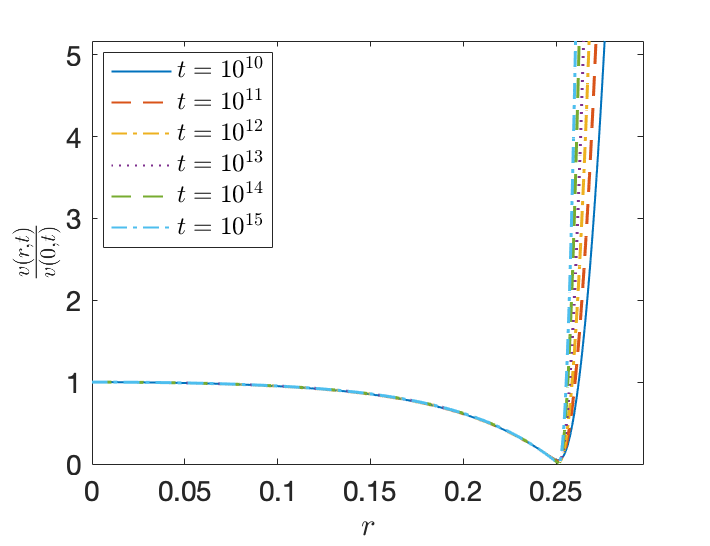}%
    \qquad%
    \includegraphics[clip=true,width=0.45\textwidth]{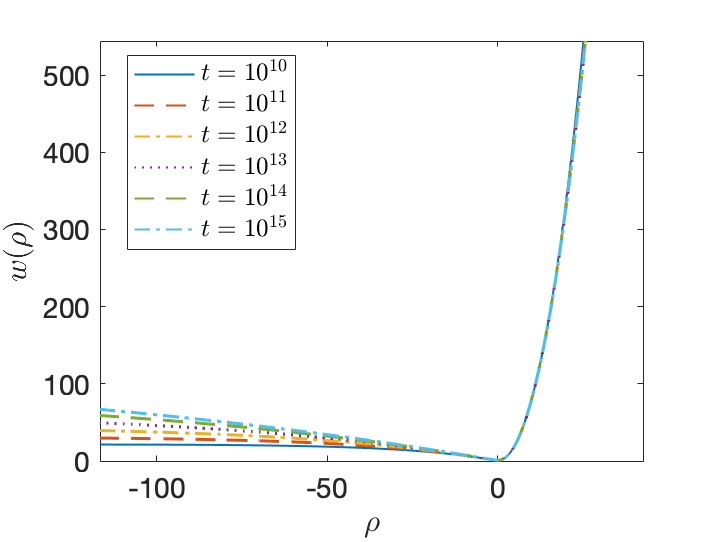}
\caption{\label{fig:crtdr} (a) Left: Central region rescaled according to $r$
vs.\ $v(r,t)/v(0,t)$ for different times (b) Right: Rescaled touchdown
region, $w$ vs $\rho$, for the same times as in (a). The
definition of $w$ and $\rho$ are given in the main text.}
\end{figure}

In the central region, we specifically make the ansatz
\begin{equation}\label{sa-central}
v(r,t)\sim t^\alpha \psi(r) 
\end{equation}
with some $\alpha<0$. The independent variable is not scaled as the region 
it spans extends from $r=0$ to near $\bar r(t)$, which is an $O(1)$ interval. 
The scaling \eqref{sa-central} can be tested by plotting $v(r,t)/v(0,t)$ in
fig.~\ref{fig:crtdr}(a), where we observe that all curves collapse near $r=0$.
Moreover, the location of the minimum of $v$, i.e. $r=\bar r(t)$, appears
to converge to a limit, which we denote by $r_*$ for future reference. Here
this limit is approximately $r_*=0.25$.

In the touchdown region, we introduce a scaling for both
variables, so that
\begin{equation}
v(r,t)\sim t^\beta \pphi(\eta), \qquad \eta:=\frac{r-r_*}{t^\gamma}, 
\end{equation}
for some $\beta$, $\gamma <0$. We test the self-similar scaling by  
scaling 
\[w:=\frac{v(r,t)}{\displaystyle \min_{r\in[0,1]} v(t)},\]
so that the minimum value of
the new function is now 1 for all $t$. 
If $\bar r(t)$ is, as before, the point where this minimum is located,
and if $r=\bar r(t)+\Delta r(t)$, with $\Delta r(t)>0$, is the point where
$w(r,t)=3$, then we define the rescaled independent variable as 
\[\rho=\frac{(r-\bar r(t))}{\Delta r(t)}.\]
Plotting $w$ vs $\rho$, we clearly see in fig.~\ref{fig:crtdr}(b) 
that around $\rho=0$, the curves collapse nicely over at least three orders of magnitude
$t=10^{10}, \ldots , 10^{15}$.

\paragraph{Similarity exponent for the central region}

\begin{figure}[p]
    \centering
        \includegraphics[width=0.29\textwidth]{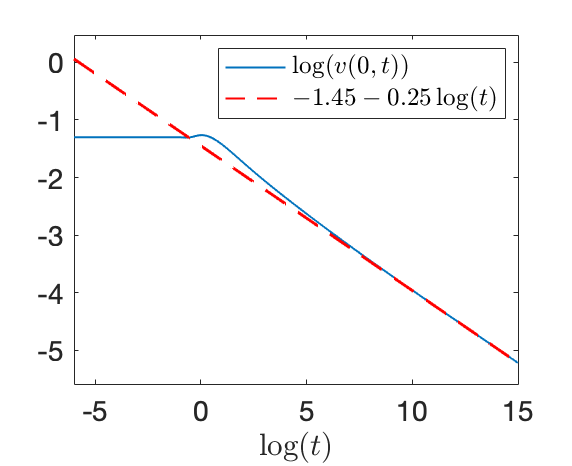}\qquad
        \includegraphics[width=0.29\textwidth]{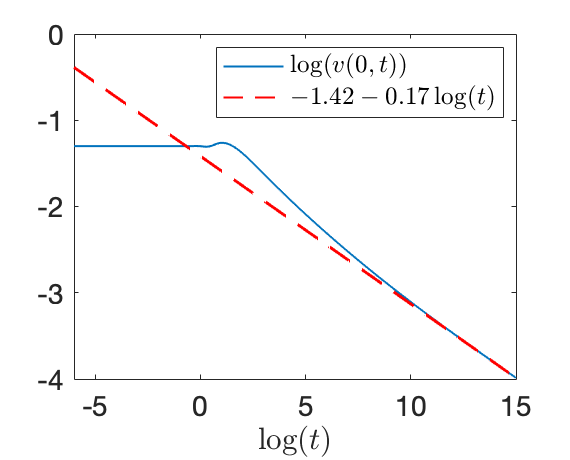}\qquad
        \includegraphics[width=0.29\textwidth]{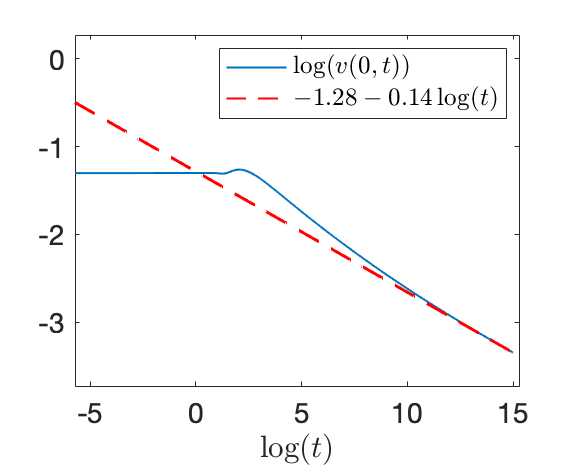}
    \caption{Linear fitting of $\log(v(0,t))$ vs $\log(t)$ for $\eps=0.1$, (a) left: $n=3$ (b), middle: $n=4$ and (c), right: $n=5$.}
    \label{fig:loglogvrzero}
\end{figure}
\begin{figure}[p]
    \centering
        \includegraphics[width=0.29\textwidth]{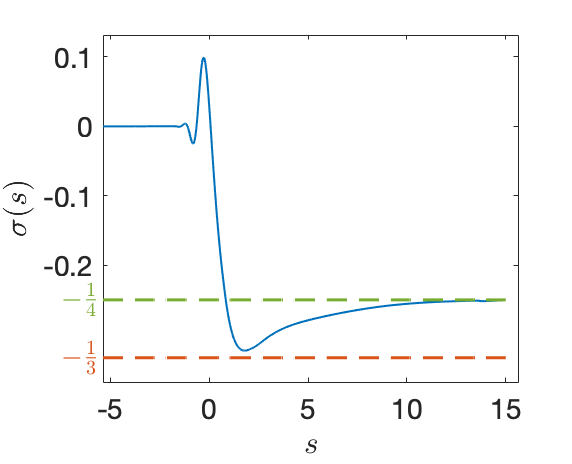}\qquad
        \includegraphics[width=0.29\textwidth]{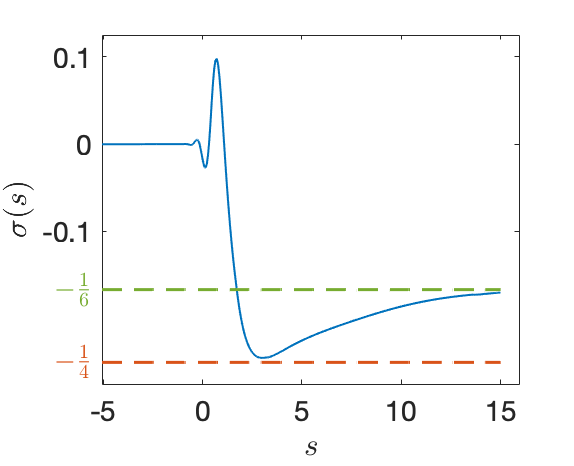}\qquad
        \includegraphics[width=0.29\textwidth]{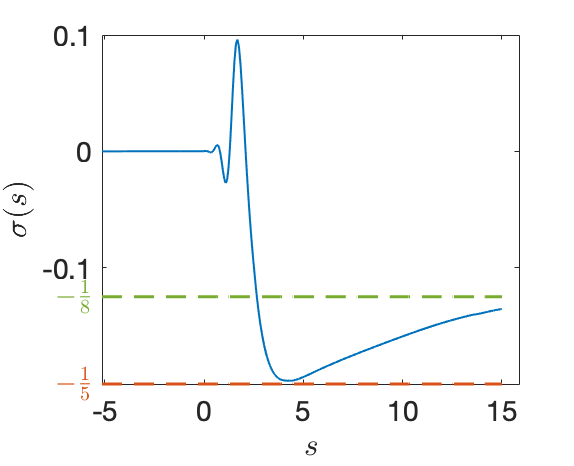}
    \caption{
Numerical approximation $\sigma$ for the power law exponent vs $s$
up to the final time $t_{final}=10^{15}$ where
we ended the simulation, and $\eps=0.1$, (a), left: for $n=3$, (b), middle:
$n=4$, (c), right: $n=5$. Definition of $\sigma$ and $s$ are given in the main text,
see \eqref{sigmaexp}.}
    \label{fig:dloglogvrzero}
\end{figure}

Next, we determine the approximate numerical values for $\alpha$, $\beta$, and
$\gamma$ from the numerical data in the central and touchdown region.
We first look at log-log plots for $v(0,t)$ for three different values of $n$.
The result is shown in fig.~\ref{fig:loglogvrzero} for $n=3,4,5$.
We see that for long times $\log v(0,t)$ is linear in
$\log t$, though it appears that two different slopes emerge at different times. 
To analyse this further, in fig.~\ref{fig:dloglogvrzero}  we
also plot 
\begin{equation}\label{sigmaexp}
\sigma(s):=\frac{\mathrm d\log v(0,t)}{\mathrm d s}, \qquad s:=\log t,
\end{equation}  
which gives the exponent of a power-law behaviour (for a pure
power law, this derivative would be exactly equal to the exponent). It is clear in fig.~\ref{fig:dloglogvrzero} 
that $\sigma$ takes a dip, with the minimum
at $-1/3$, $-1/4$, $-1/5$ for $n=3$, $4$ and $5$, respectively, suggesting 
that for general $n$, $\sigma=-1/n$ at its minimum value. However, $\sigma$ does not
remain there,
instead it increases again and then tends to $-1/4$, $-1/6$, $-1/8$ for
the three values of $n$. This indicates that for very long times, in general, we have
\[
\alpha=-\frac1{2(n-1)}.
\]

\begin{figure}[ptb]
    \centering
        \includegraphics[width=0.29\textwidth]{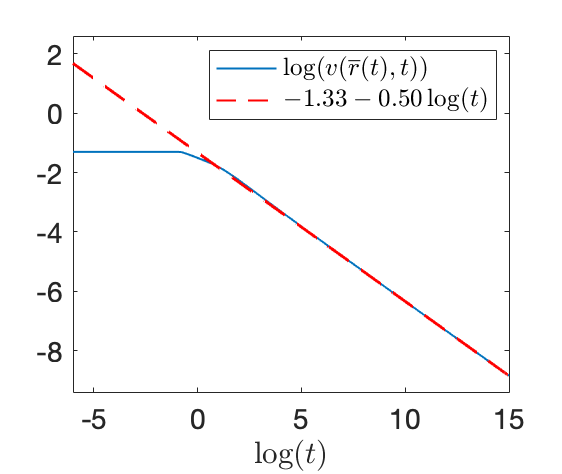}\qquad
        \includegraphics[width=0.29\textwidth]{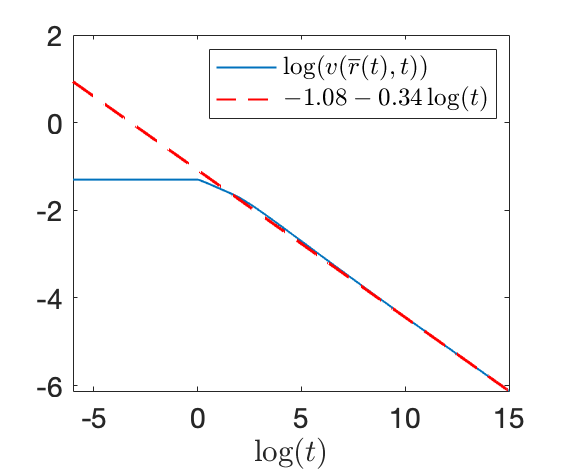}\qquad
        \includegraphics[width=0.29\textwidth]{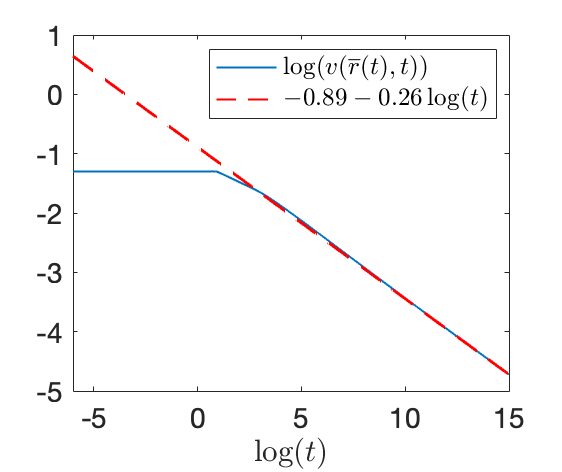}
    \caption{Linear fitting for $\log(v(\bar{r}(t),t))$ v/s $\log(t)$  for $\eps=0.1$, (a) left: $n=3$, (b) middle: $n=4$ and (c) right: $n=5$.}
    \label{fig:loglogvrstar}
\end{figure}
\begin{figure}[ptb]
    \centering
        \includegraphics[width=0.29\textwidth]{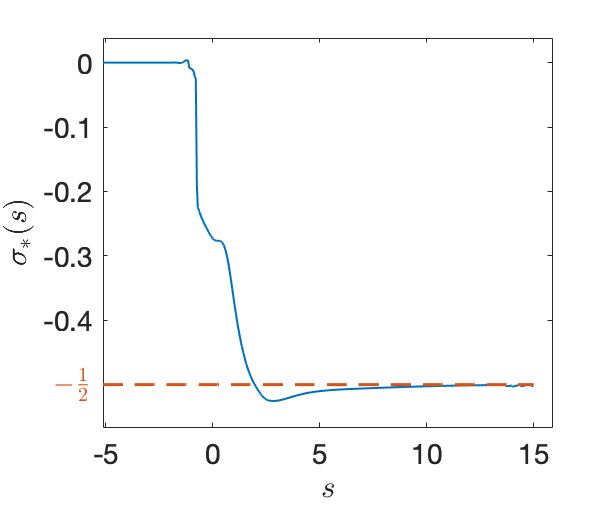}\qquad
        \includegraphics[width=0.29\textwidth]{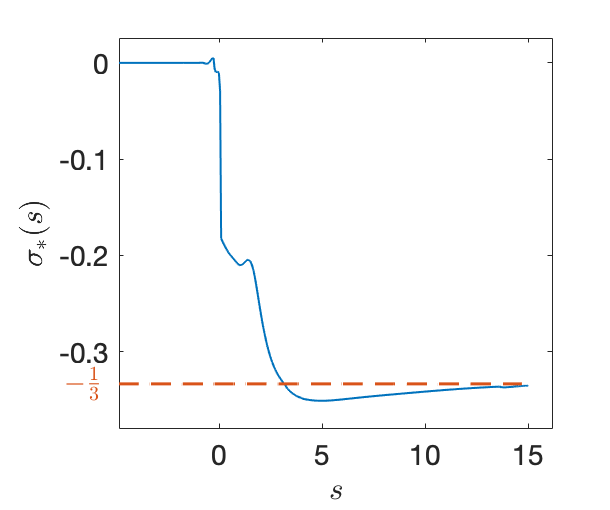}\qquad
        \includegraphics[width=0.29\textwidth]{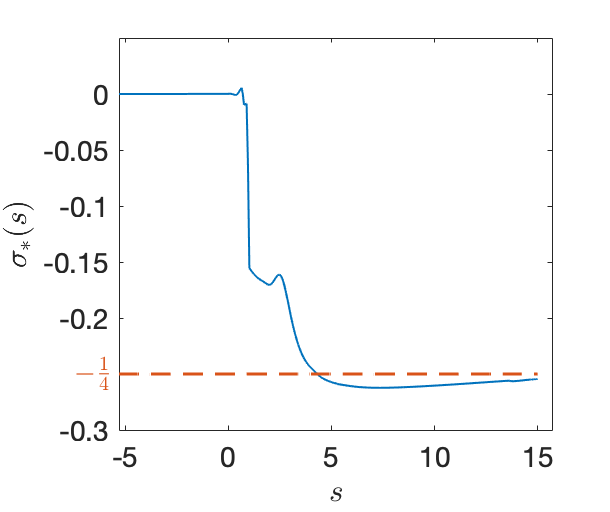}
    \caption{
Numerical approximation $\sigma_*$ for the power law exponent up to the final 
time $t_{final}=10^{15}$ where we stopped the simulation, and $\eps=0.1$, for (a), 
left: $n=3$, (b), middle: $n=4$ and (c), right: $n=5$. Definition of $\sigma_*$ and $s$ are on the main text.}
    \label{fig:dloglogvrstar}
\end{figure}

\afterpage{\clearpage}

\paragraph{Similarity exponents for the touchdown region}

In the touchdown region, the $\log$-$\log$ plot for $v(\bar{r}(t),t)$ 
in fig.~\ref{fig:loglogvrstar}
reveals that, for long
times, the evolution of this value indeed approaches a power-law
behaviour. As in the central region, approximations
for the similarity exponent $\beta$ can 
be extracted from 
\[
\sigma_*(s):=\frac{\mathrm d\log v(\bar{r}(t),t)}{\mathrm d s}.
\]
 Note that even though the touchdown region should be centered at $r_*=\lim_{t \to \infty} \bar{r}(t) $, because we do not know its value a priori, we must use instead $\bar{r}(t)$. The result is shown in fig.~\ref{fig:dloglogvrstar},
where $\sigma_*$ approaches $-1/2$, $\beta=-1/3$ and $\beta=-1/4$ for the three 
values $n=3$,
$4$, $5$, respectively. This is consistent with  
\[
\beta=-\frac{1}{n-1}.
\]
On the other hand, we obtained that $\partial_{rr} v(\bar{r}(t),t)$ tends to a constant as $t\to \infty$.
Since $\partial_{rr} v(\bar{r}(t),t)
\propto t^{\beta-2\gamma}$, this means
$\beta=2\gamma,$
and thus
\begin{equation}
\gamma=-\frac{1}{2(n-1)}.
\end{equation}


\section{Long time asymptotic analysis}\label{asymptotic analysis}

In this section we carry out a long-time asymptotic analysis of the solution
to \eqref{tp}, \eqref{mobg} via singular perturbation theory.
For this purpose, we rescale time as 
\begin{equation}\label{t2tau}
t=\tau/\delta,
\end{equation}
where $0<\delta\ll 1$ is the small parameter on which the asymptotic analysis
is based and $\tau$ is an independent time variable. 

We first formulate the leading order problem in 
each of the three regions -- central, touchdown and annular
-- solve and match them. Then we construct the composite solution and compare
it with the numerical results. 
We also note that, while $\eps$ is also small, we treat it as a fixed parameter.

For the central and touchdown regions, it is convenient to formulate the
problem in terms of the function $v=1-u$. Since we are interested in solutions
that are bounded $|u(r,t)|\leq 1$ for all $r\in (-0,1)$, $t \in (0,\infty)$, we
look for nonnegative $v$. Substituting $v$ into 
\eqref{tppde}, \eqref{tpbcat0},  \eqref{tpbcat1}, \eqref{fes}, \eqref{mobg}, 
and rescaling time according to \eqref{t2tau} we obtain
\begin{subequations}\label{eqv}
\begin{align} \label{eqv1}
    \delta \partial_\tau v & = -\frac{1}{r}\partial_r\left( r M(v) \partial_r \mu\right),\\ \label{eqv2}
    \mu&=\eps^2\left(\partial_{rr}v + \frac{1}{r} \partial_rv\right)+2(-v^3+3v^2-2v).
\end{align}
with boundary conditions 
\begin{align} \label{eqv3}
    \partial_rv(1,t) &=0,\\\label{eqv4}
    M(v(1,t))\partial_r\mu(1,t) &=0,\\\label{eqv5}
    \partial_rv(0,t) &=0,\\\label{eqv6}
    \partial_r\mu(0,t) &=0.
\end{align}
where $M(v)= v^n(2-v)_+^n$, for $n >0$.
For later uses, we record that the radial flux is given by
\begin{equation}
j=-M(v)\partial_r\mu.
\end{equation}
\end{subequations}

\subsection{Central region} \label{asymptotic outer left}

We start from \eqref{eqv} and assume, using the insight gained
from the numerical results, that $v$ can be expanded as 
\begin{equation} \label{vouter} v_{\text{central}}(r, \tau) =  
\delta^{-\alpha} \tau^{\alpha}\psi_0(r) +o(\delta^{-\alpha}),
\end{equation}
where $\alpha \in \R$ and $\psi_0$ is a nonnegative function. 
Since the solution must be bounded as $\delta \to 0$, we restrict our
attention to $\alpha<0$.

We substitute \eqref{vouter} into \eqref{eqv1}, \eqref{eqv2}, combine the equations
into a single one by eliminating $\mu$, and drop all terms that we already know are of
lower order to get 
\begin{equation} \label{vouteralpha} 
\alpha \delta^{-\alpha +1 } \tau^{\alpha-1}\psi_0  =- \delta^{-(n +1)\alpha} 
\tau^{(n+1)\alpha} \frac{2^n}{r}\partial_r \left[ r  \psi_0^n  \partial_r \left(\eps^2 
\frac{1}{r}\partial_r\left( r\partial_{r} \psi_0 \right) -4   \psi_0 \right) \right].
\end{equation}
Balancing both sides would require $\alpha=-1/n$, but this is not consistent with
the numerical results, for which $\alpha$ is clearly larger. In that case, the right
hand side of the equation dominates the left hand side, so that we obtain, after integrating
twice with respect to $r$ and 
using the Neumann boundary condition at the origin \eqref{eqv6}, 
the leading order problem 
\begin{subequations}
\begin{align} \label{finalouter11}
  \eps^2 \left( \partial_{rr}\psi_0(r) + \frac{1}{r} \partial_{r}\psi_0(r)\right)-4 \psi_0(r) &=c_1,\\\label{finalouter22}
      \partial_r\psi_0(0)&=0,
\end{align}
\end{subequations}
where $c_1$ is an unknown constant that comes from the second integration.
The general solution of \eqref{finalouter11}-\eqref{finalouter22} can be directly computed as 
\begin{equation*} 
\psi_0(r)=-\frac{c_1}{4}+ c_2 I_0\left( \frac{2}{\eps}r\right),
\end{equation*}
where $I_0$ is the \emph {modified} Bessel function of the first kind and 
$c_2$ is another unknown constant of integration.

This solution will be matched to the one in the touchdown region,
which plays the role of an inner expansion where the dependent variable
is small (in terms of $\delta\ll 1$) compared to the expansion in the central
region. Thus, $\psi_0$ must vanish at $r=r_*$, that is,
\begin{equation*}
    \psi_0(\rstar) =0,
\end{equation*}
which we use to eliminate $c_1$, giving
\begin{equation}\label{psi0c2}
\psi_0(r)=
c_2 \left(I_0\left( \frac{2}{\eps}r\right)-I_0\left( \frac{2}{\eps}\rstar\right)
\right) 
.
\end{equation}
This expression is, of course, only valid for $0\leq r\leq r_*$; 
we extend it by $\psi_0=0$ for $r_*<r\leq 1$ where this is needed (for
example in the evaluation of the composite expansion).
The remaining constant, $c_2$, represents a normalisation of $v$ that we 
keep as a parameter and that we fix when we numerically solve the problem
in the touchdown region.
For later, we record that the Taylor expansion of $\psi_0$ near $\rstar$ is
\begin{align}
    \psi_0(r)&=  a_1  (r-\rstar)
+ \order{(r-\rstar)^2},
\label{a1bessel}
\qquad
a_1 = \frac{2 c_2}{\eps}I_0'\left(\frac{2}{\eps} \rstar\right).
\end{align}


\subsection{Touchdown region} \label{asymptotic touchdown}

In this region, we introduce the independent variable 
\begin{equation} \label{eta}
    \eta = \frac{r-\rstar}{\delta^{-\gamma} \tau^{\gamma}},
\end{equation}
with $\gamma <0$, and expand
\begin{equation}\label{vtouchdown}
v_{\text{touchdown}} =  \delta^{-\beta} \tau^{\beta} \varphi_0(\eta) + o(\delta^{-\beta}),
\end{equation}
where $\beta <0$, as suggested by our previous numerical results, and $\pphi_0$ is a nonnegative function.
Dropping higher order terms and eliminating $\mu$
gives
\begin{align*}
    \delta^{1-\beta} (-\gamma \eta \tau^{-1} \partial_\eta \varphi_0 ) =& -\eps^2  \frac{\delta^{4\gamma -(n+1)\beta}\tau^{-4\gamma}  2^n}{\eta \left(\frac{\tau}{\delta}\right)^\gamma +\rstar} \partial_\eta\left( \left(\eta \left(\frac{\tau}{\delta}\right)^\gamma +\rstar\right) \varphi_0^n \partial_{\eta \eta \eta}\varphi_0\right).
\end{align*}
There are three possibilities here: either the LHS goes to zero faster than the
RHS, and hence $\beta > \frac{4\gamma -1}{n} $, the other way around, thus
$\beta <\frac{4\gamma-1}{n}$, or they balance each other with $\beta =
\frac{4\gamma -1}{n}$. From our previous numerical results,
we can infer that the first case is the relevant one.
Integrating the resulting leading order long time equation once,
we arrive at the ODE 
\begin{equation}\label{td0de}
\pphi_0^n(\eta)\partial_{\eta \eta \eta}\pphi_0(\eta)
=J,\end{equation} where the unknown flux $J$ appears as an integration constant.
This third order problem has to be matched to the central
(as $\eta\to-\infty$) and annular (as $\eta\to\infty$) regions, which both act
as larger, i.e.\ outer, layers.  
To match to a leading order
$\order1$ contribution in the annular region, the leading order term of
$\pphi_0(\eta)$ must be $\sim \eta^{\beta/\gamma}$ for $\eta\to\infty$,
so that upon scaling back into annular variables, the $\delta$-factors cancel.
From the numerical results, we have already observed that $\beta/\gamma$ is
closer to 2 than to 1, hence $\pphi_0$ grows faster than linear,
and this selects the behaviour of the annular solution to be quadratic near $r_*$.


\subsection{Annular region} \label{asymptotic outer right}

In the numerical simulations we saw that the solution evolves much slower to the right of $\overline{r}(t)$, at least compared with the central and touchdown regions. This leads us to believe that a good first approximation in the interval $(r_*,1)$ is given by the stationary problem. 
For the annular region, we expect a stationary solution to leading order, as the only
time dependence comes from the slow drainage of material from the central region.
This is supported by the numerical evidence, and therefore we let for the leading
order annular solution
\begin{equation}\label{vannular}
    v_{\text{annular}}(r,\tau) = 1- U_*(r) 
\end{equation}
Substituting into \eqref{eqv1}-\eqref{eqv6} and integrating twice, we 
obtain
\begin{subequations}\label{vouterodesys}
\begin{align}
    -\frac{\eps^2}{r}\frac{d}{dr}\left(r\frac{d U_*}{dr}\right) -2U_*(1-U_*^2) =& \mu_0, 
\label{vouterodesysa} \\
    \frac{dU_*}{dr}(1)=&0,
\label{vouterodesysb}
\end{align}
for $r\in(r_*,1)$, where $\mu_0$ is an unknown integration constant. 
The solutions can be locally expanded in a Taylor series,
which does not have a constant or linear contribution to be matchable to the 
touch down region, which grows superlinearly as observed in the previous section. Thus
\begin{equation}\label{vouterodesysc}
   U_*(r_*) = 1,\quad  \frac{dU_*}{dr}(r_*) = 0.
\end{equation}
\end{subequations}
For $r \in [0,r_*]$ we set the solution to $U_* \equiv 1$. Note that this is
exactly the problem treated by Lee et al.~\cite{lee_sharp-interface_2016}.
From \eqref{vouterodesysa} and \eqref{vouterodesysc}, 
we obtain 
the leading term in the Taylor series expansion for $v_{\text{annular}}$,
\begin{align}\label{eqn:taylor_annular}
    v_{\text{annular}} &= b_2(r-r_*)^2 +O((r-r_*)^3),\qquad
    b_2 = \frac{\mu_0}{2\eps^2}.
\end{align}


\subsection{Matching} \label{matching}
\paragraph{Central and touchdown region}

We first match the central and touchdown solutions.
First, we expand the inner expansion of $v_{\text{central}}$ as $r \to r_*$,
and rewrite the result in terms of $\eta$, giving
\begin{equation}
v_{\text{central}}=\delta^{-\alpha} \tau^\alpha a_1 (r-r_*) + \text{h.o.t.}=
\delta^{-\gamma-\alpha} \tau^{\alpha+\gamma} a_1 \eta + \text{h.o.t.},
\end{equation}
where we recall that $a_1$ 
is given in terms of $r_*$ via the modified Bessel function, 
see~\eqref{a1bessel}.
This has to be matched with
\begin{equation}
v_{\text{touchdown}}=\delta^{-\beta} \tau^{\beta} A_- \eta + \text{h.o.t.},
\end{equation}
and therefore
\begin{equation}\label{mctd}
\beta=\alpha+\gamma, \qquad A_-=a_1.
\end{equation}

\paragraph{Annular and touchdown region}

On the other side of $r_*$, we know that $v_{\text{annular}}$ has a Taylor expansion
as $r\to r_*$ that starts quadratically, hence
\begin{equation}
v_{\text{annular}}=b_2 \delta^{-2\gamma} \tau^{2\gamma} \eta^2 + \text{h.o.t.},
\end{equation}
where again $b_2$ is known from \eqref{vouterodesys} in terms of $r_*$.
Thus, the expansion of the solution in the touchdown region
at $\eta\to\infty$ must also be quadratic
\begin{equation}
v_{\text{touchdown}}=\delta^{-\beta} \tau^{\beta} A_+ \eta^2 + \text{h.o.t.}
\end{equation} 
and completing the matching requires
\begin{equation}\label{mtda}
\beta=2\gamma, \qquad A_+=b_2.
\end{equation}

\paragraph{Matching of the flux between central and touchdown region}
So far, we have only got two relations for $\alpha$, $\beta$ and $\gamma$,
we need one more to completely fix the similarity exponents. 
One can obtain a partial mass conservation condition in the interval $(0,r_*)$ by multiplying \eqref{eqv1} by $r$, then integrating in $r \in (0,r_*)$ and using the boundary condition at $r=0$ \eqref{eqv6}. This gives 
\begin{equation*}
    \int_0^{r_*} \partial_t v(r,t) r dr = - r_* M(v(r_*,t)) \partial_r \mu (r_*,t). 
\end{equation*}
This means that the rate of change in mass in the interval $(0,r_*)$ is
equivalent to the flux at $r_*$. 
We rescale the right hand side into touch-down variables and use the leading
order asymptotic solutions \eqref{vouter} and \eqref{vtouchdown}, to obtain 
\begin{equation*}
    -\alpha \delta^{-\alpha +1} \tau^{\alpha -1} \frac{c_2 r^{*^2}}{2} I_2\left(\frac{2r_*}{\varepsilon}\right)= - r_*  \delta^{-\beta(n+1) + 3 \gamma} \tau^{\beta(n+1) - 3 \gamma} 2^n \varepsilon^2 J,
\end{equation*}
and, therefore, matching requires 
\begin{align}
    \label{mj}
    \alpha - 1 &= \beta (n+1) - 3 \gamma,\\\label{Jc2}
    J&=\alpha  \frac{c_2 r_*}{2^{n+1} \varepsilon^2 } I_2\left(\frac{2r_*}{\varepsilon}\right).
\end{align}
We note that this could also be obtained from matching at higher order in the expansions instead of using the mass conservation derived from the equation.
The solution to \eqref{mctd}, \eqref{mtda}, \eqref{mj} is
\begin{equation}
\alpha=\gamma=-\frac{1}{2(n-1)}, \quad
\beta=-\frac{1}{n-1}.
\end{equation}
These are exactly the values that we observed in the numerical
results in section.


\subsection{Solution in the touchdown region}
We analyse the touchdown region in more detail to ensure that a
solution can be obtained, at least numerically, 
that satisfies all the matching conditions. In particular, we 
carry out an overall degree of freedom count, and briefly explain
how we solve for $\phi_0$.

The expansions of solutions $\phi_0$ 
of \eqref{td0de}
for large negative and positive arguments can be obtained from
the literature, see e.g.\ the systematic study of such expansions
for thin-film type equations in \cite{boatto_traveling-wave_1993},
but to be self-contained, we give a derivation in the appendix for
the case of linear leading order as $\eta\to-\infty$ and quadratic leading
order for $\eta\to\infty$. In summary, we have
\begin{subequations}
\begin{align}
\label{finaltouchdown2}
    \pphi_0(\eta) &= \begin{cases} A_- \eta + \frac{JA_-^{-n}(-\eta)^{3-n}}{(n-1)(n-2)(n-3)} +B_- +O(1) & \text{ if } 2<n < 3,\\
    A_-\eta +\frac{J}{2A_-^3} \ln(-\eta) + B_- + O(1) & \text{ if } n=3,\\
    A_- \eta +B_- + O(1)  & \text{ if } n > 3,\end{cases}  ,
\intertext{as $\eta \to -\infty$, and}
\label{finaltouchdown3}
    \pphi_0(\eta) &= A_+ \eta^2 +B_+\eta + C_+ + O(1)
\quad \text{ for } n>3/2  \text{\quad as } \eta \to \infty.
\end{align}
\end{subequations}
where $A_\pm$, $B_\pm$, $C_+$, $J$ are unknown constants. The limitations on $n$ arise
from the requirement that the correction terms must be asymptotically small compared
to the leading order term. In particular the expansions for large negative arguments
is only valid for $n>2$ and sets the lower bound for the
$n$ that we consider in this study. 

The degree of freedom count is as follows: We have 6 unknown constants,
$A_\pm$, $B_\pm$, $C_+$ and $J$. Two of them ($A_+$ and $A_-$) have been fixed by
matching to the appropriate outer problems, though this introduces additional unknowns,
which we will return to later.  The
third order differential equation takes care of another three degrees of
freedom, so that only one is left. This degree of freedom is the result of
the ODE \eqref{td0de} being autonomous, and represents an arbitrary shift of
the solution. This shift is fixed by the requirement that $\eta$ is defined in
\eqref{uinit} through the position $r_*$ where $v$ touches down (for $t\to\infty$).
Hence $\varphi_0$ must have its minimum at $\eta=0$.

Returning first to $A_-$, we see that matching specifies this constant in
terms of $c_2$ and $r_*$ see \eqref{a1bessel} and \eqref{mctd}, but
we can use \eqref{Jc2} to eliminate $c_2$. 
For $A_+$, the equations \eqref{mtda} and \eqref{eqn:taylor_annular} introduce
a dependence on $\mu_0$. However, both $r_*$ and $\mu_0$ are completely determined
by solving the leading order problem in the annular region
\eqref{eqn:taylor_annular}, \eqref{vouterodesys} if this is supplemented
by an overall mass constraint. Hence after matching all layers, the solutions
in all regions are completely determined.

The solution strategy for the touchdown region is as follows. We consider 
\eqref{td0de} first with the condition \eqref{eqn:taylor_annular}.
Since $B_-$ represents a translation in $\eta$ we do not need to enforce its
value and only impose the leading order behaviour $A_- \eta$. By rescaling
\[\pphi_0(\eta) = c \phi_0(y), \quad \eta = d y,\]
with
\[c := \left(\frac{J}{(-A_-)^3}\right)^{\frac{1}{n-2}}>0, \quad d := \left(\frac{J}{(-A_-)^{n+1}}\right)^{\frac{1}{n-2}}>0.\]
we see that $\phi_0$ satisfies the following parameter free problem
\begin{subequations}\label{phi0norm}
\begin{align}
    \phi_0^n(y)\partial_{yyy}\phi_0(y) &=1, \quad y \in (-\infty, \infty), \\
    \phi_0(y) &=  -y \text{ as } \eta \to -\infty .
\end{align}
\end{subequations}
and read off $\kappa = \partial_{yy}\phi_0(y_0)$ from the numerical solution. 
Scaling
back and using \eqref{mtda}, \eqref{eqn:taylor_annular} gives 
\[2\frac{\mu_0}{\eps^2}=2A_+=\kappa \frac{c}{d^2}.\]
This is essentially an equation between $J$, $A_-$ and $\mu_0$, the latter
being fixed by the annular region. Replacing $A_-=a_1$ by the second equation
in \eqref{a1bessel} and using \eqref{Jc2},
we obtain the following expression for $c_2$
\[
c_2 = \left(\frac{\mu_0^{n-2} \eps I_2\left(\frac{2r_*}{\eps}\right) }{2^{3n+1}(n-1)  \kappa^{n-2} I_1\left(\frac{2r_*}{\eps}\right)^{2n-1}}\right)^{\frac{1}{2(n-1)}},
\]
which now fixes $c_2$ once $\mu_0$ and $r_*$ have been obtained by solving the
annular problem. 
Numerically, \eqref{phi0norm} is solved on a large truncated domain and then extended,
where necessary, to an infinite domain by using the expansions 
\eqref{finaltouchdown2}, \eqref{finaltouchdown3}.

\subsection{Solution in the annular region}

The solution in the annular region also needs to be obtained numerically. As
noted in \cite{lee_sharp-interface_2016} we need an extra condition to solve
for the unknown $\mu_0$, which may be obtained by either fixing the position of the
interface or adding a mass constraint. We choose the latter and impose
\begin{equation}\label{massconstraint}
\int_{r_*}^r U_*(s)s \,ds=m_0-\frac{r_*^2}{2},
\end{equation}
with $m_0$ is the initial mass, as defined in \eqref{mzero}.
The problem \eqref{vouterodesys}, \eqref{massconstraint} 
then solved by picking $r_*$, solving all conditions except for \eqref{massconstraint}
using the Matlab solver bvp4c, and then iterating over $r_*$ until the mass constraint
is satisfied, too. Practical details, such as the conversion into a boundary value problem
for a system of first order ODEs, are discussed in an appendix. 
For the problem here, we obtain a solution with the required mass for $r_*= 0.2516$.

\subsection{Composite approximation} \label{composite}

\begin{figure}[tb]
    \centering
    \includegraphics[clip=true,width=0.45\textwidth,height=0.36\textwidth]{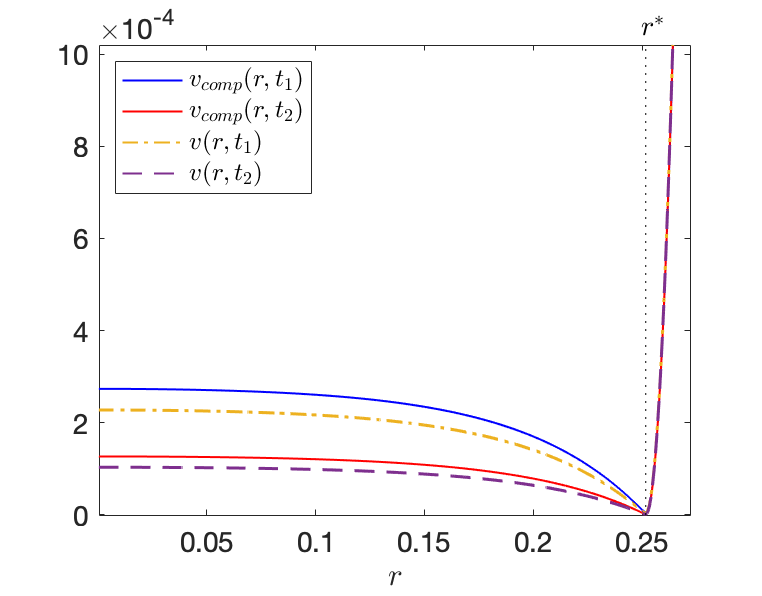}
    \includegraphics[clip=true,width=0.45\textwidth,height=0.36\textwidth]{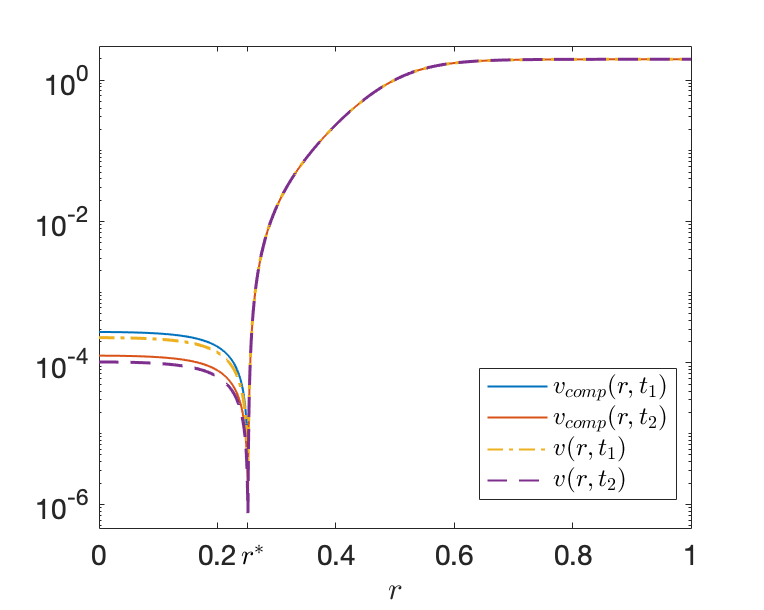}
    \caption{(a) Left: Plot of $v_{\text{comp}}$ v/s $r$ compared to numerical PDE solution
$v$, for $t_1=10^{13}$, $t_2=10^{15}$, $\Delta r=10^{-6}$, $n=4$, $\eps=0.1$ and $r_*= 0.2516$.
(b) Right: Comparison of $v_{\text{comp}}$ with the PDE solution, on a scale that 
shows the complete solution. Parameter values and times carry over
from (a).
}
    \label{fig:compsol1}
\end{figure}

We now construct the composite approximation from the asymptotic solutions
found in the previous sections. We add each of the approximations - central, touchdown and annular - in the same variables and subtract common terms, in other words
\begin{align}
v_{\text{comp}}(r,t)=& v_{\text{central}}(r,t) + v_{\text{touchdown}}(r,t)
+v_{\text{annular}}(r,t)
\notag\\ & \quad - A_- t^{\alpha} (r-\rstar)_-
- A_+ (r-\rstar)_+^2 .\label{comp3}
\end{align}
The subscript ``$-$'' (``$+$'') means that we take the value of the brackets
where it is negative (positive), and zero elsewhere.  

In Figures \ref{fig:compsol1}(a) and (b) we present the leading order of the composite
expansion at two times, $t_1=10^{13}$ and $t_2=10^{15}$, and compare them with the
numerical PDE solution for $r\in [0,1]$ at both times. We use a semi-log plot in \ref{fig:compsol1}(b) so that all relevant parts of the solution, which differ
by several orders of magnitude, can be shown in a single plot. Note that we only use the PDE solution for comparison and the only external data that is fed into the 
computation of composite expansion is the mass of the initial condition $m_0$. 
The composite solution closely follows the numerical solution of the full
problem for $v$ in the whole interval $[0,1]$. In particular, around the
minimum of the solution and in the annular region the agreement is excellent,
as shown in Fig.~\ref{fig:compsol1}(a) and (b), at times $t=t_1$ and $t_2$.
In fact, the overall agreement becomes better at the later time, $t=t_2$, as
the absolute error
\[
\text{error}(t)= \max_{r\in [0,1]}|v(r,t)-v_{\text{comp}}(r,t)|
\]
decreases by a factor of 7.4 between $t_1$ and $t_2$.

\section{Discussion and outlook}

While the spontaneous appearance of touchdowns, that is, the emergence of points
$\{ M(u)=0 \}$ for the degenerate Cahn-Hilliard equation is a well-known
phenomenon in the development of numerical algorithms for these types of PDEs,
so far, they have been addressed by ad-hoc measures such as adding a small
positive constant to the mobility, see for example \cite{torabi_new_2009}.
In this study, we systematically investigate how such points of vanishing mobility arise.
Important results of this work are that small energy solutions to the
degenerate Cahn-Hilliard problem \eqref{tp} touch down, either in finite or
infinite time, and that for $n>2$ and practically relevant
tanh-like initial data, the solution converges to a long-time profile that
touches down in infinite time, suggesting this is the generic behaviour for
large enough values of $n$.

Our asymptotic analysis also revealed that $n=2$ itself is a critical case and
the values of $n\leq 2$ require a separate investigation, which will appear in
an upcoming paper. There it is shown that asymptotic approximations for
infinite-time touchdown solutions can be found also for $1/2<n\leq 2$, but
they become increasingly fragile. This suggests that for $n$ closer to $1/2$,
finite-time touchdown prevails in the degenerate Cahn-Hilliard problem
\eqref{tp}.

Future work could investigate fully 2D situations and non-convex interfaces
between the phases, and the impact of touchdown regions on the sharp interface
evolution. For finite time touchdown, once the solution has reached $|u|=1$,
further questions arise about how to continue the solution, noting
that already more than one weak solution concept has been considered in the
literature \cite{elliott_cahnhilliard_1996, dai_weak_2015}.  Also, thin film
theory could provide insight into how to construct numerical schemes that
maintain $|u|<1$, e.g. \cite{zhornitskaya_positivity-preserving_2000,
grun_nonnegativity_2000}. 

More generally, our study underscores the close connection of degenerate
Cahn-Hilliard equations with thin film models, and the still untapped potential
for using not only the rigorous proofs, but also the asymptotic results from
thin film theory to understand singularity formation in degenerate
Cahn-Hilliard models. Inside a spherically symmetric phase domain (more
generally, convex domain) in the late stages of phase separation, the
degenerate Cahn-Hilliard model is closely related to the thin film equation
with pressure boundary conditions \cite{bertozzi_singularities_1994} and hence
the points with vanishing mobility can be analysed in a similar way. This opens
up a rich source of analytical tools, both asymptotic and rigorous, that we can
apply to degenerate Cahn-Hilliard problems \cite{constantin_droplet_1993,
constantin_singularity_2018, beretta_nonnegative_1995, bertozzi_symmetric_1996}.


\textbf{Funding:} This work was supported by the EPSRC Centre for Doctoral Training
in Partial Differential Equations: Analysis and Applications,
through EPSRC grant EP/L015811/1. CP acknowledges support from the Chilean National
Agency for Research and Development (ANID) through the  Scholarship Program
DOCTORADO BECAS CHILE/2018 - 72190029.\\

\textbf{Acknowledgements:} CP and AM thank Amy Novick-Cohen and Peter Howell for useful discussions.


\begin{appendix}

\section{Expansion of $\phi_0$ for large arguments}

\paragraph{Expansion at $\eta\to-\infty$}

We let $x= -\eta$, and look for an expansion of $\phi(x):=\pphi_0(\eta)$ 
such that 
\begin{equation}\label{phi}
    \phi^n \phi''' =-J,
\end{equation}
as $x\to \infty$, where $'$ denotes derivatives with respect to $x$.

Starting with the expansion $\phi(x)=Ax+\xi(x)$ with $\xi\ll x$ as $x \to \infty$ gives,
upon substituting this ansatz into \eqref{phi} and balancing, the following corrections
\begin{equation}\label{xigen1}
\xi(x) = \frac{-J}{A^n(n-1)(n-2)(n-3)}x^{3-n}+B,
\end{equation}
where $B$ is a constant, provided $n\neq 3$. 
We note that the first term dominates the second if $n<3$, and vice versa if $n>3$.
Consistency between the leading order and correction requires that $3-n<1$,  i.e. $n>2$,
because only then we have that $\xi\ll x$. Moreover, if $n=3$, the confluence of $x^{3-n}$ and
the contribution from $x^p$ with $p=0$ produces a logarithmic term, that is,
\begin{equation}\label{xilog}
\xi (x)= \frac{-J}{2A^3} \ln(x) +B.
\end{equation}
Returning to the original variables then gives \eqref{finaltouchdown2}.

\paragraph{Expansion at $\eta\to\infty$}
For $\eta \to \infty$ we are looking for an expansion of solutions
of \eqref{td0de} 
starting with a quadratic term, where $'$ denotes derivatives with respect to $\eta$. We make the ansatz
\begin{equation*}
    \pphi_0 (\eta)= A\eta^2 +\xi(\eta),
\end{equation*}
with $\xi\ll \eta^2$ and $A$ a constant. 
Introducing this into the differential equation, we obtain
\begin{equation}\label{tildephi}
\xi (\eta) = D\eta+ \frac{-J}{A^n(2n-1)(2n-2)(2n-3)}\eta^{3-2n} +E,
\end{equation}
where $D$ and $E$ are unknown constants.  For  $n\leq 1/2$, the second term grows faster
than quadratic and therefore the expansion of $\phi_0$ is not consistent.
Hence we require $n>1/2$.  Note that the order of the terms in \eqref{tildephi} changes as the values
$n=3/2$ and $n=1$. If $n>3/2$, the term proportional to $J$ decays for large $\eta$ and
therefore in this case the first three terms in the $\eta\to\infty$ expansion of $\phi_0$ are
as claimed in \eqref{finaltouchdown3}.

\section{Details of the numerical approach in the annular region}

We define 
\begin{equation*}
    S(r):= \int_{r_*}^r U_*(s)s \,ds.
\end{equation*}
Then the system \eqref{vouterodesys} with \eqref{massconstraint} can be written as
\begin{subequations}\label{ann}
\begin{align}\label{annular1}
     S' =& U_* r,\\\label{annular2}
     U_*' =& W,\\\label{annular3}
    W' =& \frac{2}{\eps^2} U_*(U_*^2-1) - \frac{\mu_0}{\eps^2} - \frac{1}{r} W,\\\label{annular4}
    \mu_0' =& 0,\\\label{annular5}
   U_*(r_*) =& 1, \quad   W(r_*) = 0, \quad  S(r_*) =0\\\label{annular6}
   W(1) =& 0, \quad S(1) =  m_0-\frac{r_*^2}{2}.
\end{align}
\end{subequations}

For an initial choice of $r_* \in (0,1)$ we solve \eqref{ann}
except for the condition on $S(1)$
using bvp4c from Matlab, and then iterate over $r_*$ until the condition on $S(1)$
is satisfied, too.

As initial guess for bvp4c we use
\begin{align*}
    S(r)= -\int_0^1\tanh \left(\frac{r-0.5}{\eps}\right) r dr,& \quad U_*(r)= -\tanh \left(\frac{r-0.5}{\eps}\right),\\
    W(r)=  -\frac{1}{\eps}\text{sech}\left(\frac{r-0.5}{\eps}\right)^2,& \quad \mu_0(r) = 1,
\end{align*}
which proved to be sufficient so that the algorithm converges.

\end{appendix}

\bibliographystyle{abbrv}
\bibliography{ref}

\end{document}